\documentclass{gtart_a}
\pdfoutput=1

\usepackage{pinlabel}
\usepackage[all]{xy}

\title[Homogeneous coordinate rings and mirror symmetry]{Homogeneous
coordinate rings and\\mirror symmetry for toric varieties}
\author[M Abouzaid]{Mohammed Abouzaid} 
\givenname{Mohammed}
\surname{Abouzaid}
\address{Department of Mathematics\\University of Chicago\\\newline
 Chicago, IL 60637\\USA}
\email{mabouzai@math.uchicago.edu}
\urladdr{}

\volumenumber{10}
\issuenumber{}
\publicationyear{2006}
\papernumber{26}
\lognumber{0683}
\startpage{1097}
\endpage{1156}

\doi{}
\MR{}
\Zbl{}


\received{26 November 2005}
\revised{3 May 2006}
\accepted{1 June 2006}
\published{}
\publishedonline{24 August 2006}
\proposed{Jim Bryan}
\seconded{Lothar G\"ottsche, Simon Donaldson}
\corresponding{}
\editor{}
\version{}

\arxivreference{math.SG/0511644}

\let\xysavmatrix\xymatrix
\def\xymatrix{\disablesubscriptcorrection\xysavmatrix}
\AtBeginDocument{\let\bar\wbar\let\tilde\wtilde\let\hat\what}
\def\SetFigFont#1#2#3#4#5{\small}
\def\acite#1*#2{\cite[#2]{#1}}
\def\adm{\mathrm{adm}}
\let\mathscr\mathcal

\numberwithin{equation}{section}

\newcommand{\im}{\operatorname{Im}}

\newcommand{\into}{\hookrightarrow}

\newcommand{\bC}{{\mathbb C}}

\newcommand{\bP}{{\mathbb P}}

\newcommand{\bR}{{\mathbb R}}

\newcommand{\bT}{{\mathbb T}}

\newcommand{\bZ}{{\mathbb Z}}

\newcommand{\cA}{{\mathcal A}}

\newcommand{\cC}{{\mathcal C}}

\newcommand{\cH}{{\mathcal H}}

\newcommand{\cL}{{\mathcal L}}
\newcommand{\cM}{{\mathcal M}}

\newcommand{\cO}{{\mathcal O}}

\newcommand{\sL}{{\mathscr L}}

\newcommand{\dd}{\partial}

\newcommand{\Hom}{\operatorname{Hom}}

\newcommand{\Log}{\operatorname{Log}}
\newcommand{\Proj}{\operatorname{Proj}}

\newcommand{\grad}{\operatorname{grad}}

\newcommand{\ctorus}[1]{(\mathbb{C}^{\star})^{#1}}

\makeatletter
\def\cnewtheorem#1[#2]#3{\newtheorem{#1}{#3}[section]
\expandafter\let\csname c@#1\endcsname\c@thm}

\newtheorem{thm}{Theorem}[section]
\cnewtheorem{cor}[thm]{Corollary}
\cnewtheorem{lem}[thm]{Lemma}
\cnewtheorem{prop}[thm]{Proposition}

\theoremstyle{remark}
\cnewtheorem{defin}[thm]{Definition}
\cnewtheorem{rem}[thm]{Remark}
\cnewtheorem{rems}[thm]{Remarks}
\newtheorem*{claim}{Claim}

\makeatother  

\newcommand{\noproof}{\hfill \qedsymbol}

\begin{document}

\begin{asciiabstract}
Given a smooth toric variety X and an ample line bundle O(1), we construct
a sequence of Lagrangian submanifolds of (C^*)^n with boundary on a level
set of the Landau-Ginzburg mirror of X. The corresponding Floer homology
groups form a graded algebra under the cup product which is canonically
isomorphic to the homogeneous coordinate ring of X.
\end{asciiabstract}

\begin{htmlabstract}
Given a smooth toric variety X and an ample line bundle
O(1), we construct a sequence of Lagrangian submanifolds of
(<b>C</b><sup>*</sup>)<sup>n</sup> with boundary on a level set of
the Landau&ndash;Ginzburg mirror of X. The corresponding Floer homology
groups form a graded algebra under the cup product which is canonically
isomorphic to the homogeneous coordinate ring of X.
\end{htmlabstract}

\begin{abstract}
Given a smooth toric variety $X$ and an ample line bundle
$\mathcal{O}(1)$, we construct a sequence of Lagrangian submanifolds
of $(\mathbb{C}^{\star})^n$ with boundary on a level set of the
Landau--Ginzburg mirror of $X$. The corresponding Floer homology
groups form a graded algebra under the cup product which is canonically
isomorphic to the homogeneous coordinate ring of $X$.
\end{abstract}

\maketitle

\section{Introduction}

In this paper we give some evidence for M Kontsevich's homological
mirror symmetry conjecture \cite{kontsevich} in the context of toric
varieties.  Recall that a smooth complete toric variety is given by a
simplicial rational polyhedral fan $\Delta$ such that $|\Delta| =
\bR^{n}$ and all maximal cones are non-singular (Fulton
\acite{fulton}*{Section 2.1}).  The convex hull of the primitive
vertices of the $1$--cones of $\Delta$ is a convex polytope which we
denote by $P$, containing the origin as an interior point, and may be
thought of as the Newton polytope of a Laurent polynomial $W \co
\ctorus{n} \to \bC$. This Laurent polynomial is the Landau--Ginzburg
mirror of $X$.

Our construction will rely on choosing an ample line bundle $\cO(1)$
on $X$.  On the complex side, we know that sections of this line
bundle are given by lattice points of a polytope $Q$ determined by
$\cO(1)$.  Further, we can recover the original toric variety as
\begin{equation}\label{proj}  X = \Proj \left( \bigoplus_{j \geq 0} \bC^{j Q \cap \bZ^{n}} \right), \end{equation} with product given by linearly extending the formula
\[ x_{k} \otimes y_{l} \to x_{k} + y_{l} .\]
For the bulk of the paper, it will be more convenient to identify the
integral points of $jQ$ with the $\frac{1}{j}$ integral points of $Q$,
and re-write the product as
\begin{equation} \label{fam}
x_{k} \otimes y_{l} \to \frac{k x_{k} + l y_{l} }{k+l}.
\end{equation}
On the symplectic side, we study the Floer homology of Lagrangian
submanifolds of $\ctorus{n}$ with boundary on $W^{-1}(0)$.  $\cO(1)$
will induce a subdivision of $P$, and hence a tropical degeneration of
$W^{-1}(0)$ in the sense of Mikhalkin.  In order to simplify our computations,
we will have to replace $\cO(1)$ by a sufficiently high power (see
\fullref{int}).  The combinatorial data which determines $\cO(1)$
will then allow us to construct Lagrangians which we suggestively call
$\cL(j)$ for every integer $j$.  We will write $\cL$ for $\cL(0)$. Our main result is the following:
\begin{thm} \label{main}
The Floer cohomology groups $HF^{0}(\cL,
\cL(j))$ for $j \geq 0$ form an associative algebra under the cup product which is canonically isomorphic
to the homogeneous coordinate ring of $X$ equipped with the line
bundle $\cO(1)$.
\end{thm}

We now summarize the contents of this paper.  In Section 2, we will
introduce the notion of admissible Lagrangian (originally due to
Kontsevich), and explain why their Floer theory is well defined.  In
Section 3, we will review the results of Mikhalkin \cite{mikhalkin}
which we will need, and explain how the datum of a line bundle on a
toric variety gives rise to a tropical degeneration of its mirror.
Ignoring boundary problems, we will also introduce some flat
Lagrangian submanifolds of $T^{*} \bR^{n} / \bZ^{n} $ whose Floer
cohomology groups correspond to holomorphic sections of the line
bundles $\cO(j)$.  In Section 4, we will give an explicit construction
of a smooth symplectic submanifold which interpolates between the
complex hypersurface and its ``tropical'' counterpart, as well as a
construction of the Lagrangian submanifold $\cL$ which corresponds to
the structure sheaf of $X$.  In Section 5, we will use the Lagrangians
introduced in Section 3 to construct the admissible Lagrangians
$\cL(j)$ that appear in the main theorem, and complete its proof.
There is a minor change of notation which occurs in the middle of the
paper as explained in \fullref{switch}.

We would like to comment on some related work.  We do not discuss any homological algebra in this paper, but the result we prove establishes the existence of a functor from the Fukaya category of the mirror of $X_{\Delta}$ to the category of coherent sheaves on $X_{\Delta}$.  Let us briefly sketch the construction.  The functor takes every Lagrangian to the graded vector space $\bigoplus_{j \geq 0} HF^{0}(\Lambda, \cL(j))$.  By \fullref{main} this is naturally a graded module over the homogeneous coordinate ring of $X$.  A classical result of Serre implies that the category of coherent sheaves is a quotient of the category of graded modules.  

In order to prove that this functor yields an equivalence of derived categories, one would need to understand how the $A_{\infty}$ structures on both sides affect the construction.  For $\bC \bP^{2}$, this can be done explicitly using essentially the methods of this paper and Beilinson's description of the category of coherent sheaves on projective spaces \cite{beilinson}, and we expect the computation to extend to $\bC \bP^{n}$.  However, for a general toric variety, it is not clear how to obtain such an explicit description of the category of coherent sheaves.  In an upcoming paper, we will establish such an equivalence of categories, but the methods are unfortunately much less explicit than those used here, and pass through Morse theory.  Also, if we consider the Lagrangians $\partial L(j) \subset M$, we can use the ideas of this paper to give evidence for the mirror conjecture, this time in the case of Calabi--Yau hypersurfaces in toric varieties.

\subsection*{Acknowledgments}
I would first like to thank my advisor Paul Seidel  for originally suggesting the use of tropical geometry as an approach to Kontsevich's conjecture, and for helpful discussions throughout my work on this project.  I also thank Gabriel Kerr for teaching me many things about toric varieties, Kevin Costello for useful comments on an early draft, and Ivan Smith for an important conversation about this work while it was still in maturation.  I finally want to acknowledge the referee whose careful reading of an earlier version unearthed many errors, not all of them minor, and whose suggestions for addressing the more significant issues have been particularly helpful.

\section{Floer (co)-homology for Lagrangians with boundary along a complex hypersurface}

\subsection{Admissible Lagrangians}
We begin by observing that many properties of closed Lagrangian
manifolds can be extended to the case where the Lagrangians have
boundary.  For example, the cotangent bundle of any manifold with
boundary is itself a manifold with boundary which carries the usual
symplectic structure.

Recall that the restriction of the cotangent bundle to the boundary
carries a canonical oriented rank--$1$ trivial sub-bundle $E_{\partial
L}$ which annihilates the tangent space to the boundary and induces
the appropriate co-orientation on $\partial L$.  Note that we can think of this
sub-bundle as lying in the restriction  of the tangent space of $T^{*} L$ to the boundary of the zero section.  The proof of Weinstein's
neighbourhood theorem extends to this setting to show that Lagrangian
submanifolds with boundary, with a choice of a rank--$1$ oriented
sub-bundle at the boundary, are locally modeled after the cotangent
bundle with its canonical sub-bundle at the boundary.
\begin{lem} \label{weinst}
Let $L$ be a Lagrangian submanifold of a symplectic manifold $N$, and
let $E$ be an oriented rank--$1$ sub-bundle of the symplectic orthogonal complement of
$T \partial L$ such that the pairing
\[ TL |_{\partial L} \otimes E \to \bR \]
induced by the symplectic form is non-degenerate and yields the
appropriate co-orientation on $\partial L$.

Inside a sufficiently small neighbourhood of $L$ in $N$, there exists
a full dimensional submanifold with boundary $(V_{L},\partial V_{L})$,
such that the inclusion $(L, \partial L) \subset (V_L, \partial V_L)$ satisfies the following properties:
\begin{itemize}
\item The restriction of $T \partial V_L$ to $\partial L$ contains the
sub-bundle $E$.
\item There exists a symplectomorphism $(V_L, \partial V_L) \to
(T^{*}L, T^{*}L|_{\partial L})$ identifying $L$ with the zero section of its cotangent bundle and $E$ with the canonical sub-bundle $E_{\partial L}$ and a projection to $(L, \partial L)$ such that the following diagram commutes:
$$
\xymatrix{ (V_L, \partial V_L) \ar[dr] \ar[rr] & & (T^{*}L,
 T^{*}L|_{\partial L}) \ar[dl] \\ &(L, \partial L) & }
$$\qed
\end{itemize} 
\end{lem}

As usual, this allows us to reduce problems about the
topology of nearby Lagrangian submanifolds to questions about closed
forms on $L$.  We will find the following lemma particularly useful:
\begin{lem} \label{C1-close} Let $L'$ and $L$ be two Lagrangian submanifolds of $N$ which have the
same boundary.  Let $V_L$ be a submanifold of $N$ which satisfies the
conditions of \fullref{weinst}.  If $L'$ is transverse to  $\partial V_L$ and there is a neighbourhood of $\partial L'$ in $L'$ which is
contained in $V_L$, then there exists a Lagrangian submanifold $L''$
which satisfies the following conditions:
\begin{itemize}
\item $L''$ is Hamiltonian isotopic to $L$.
\item $L''$ agrees with $L'$ in a sufficiently small neighbourhood of
$\partial L'$.
\item $L''$ agrees with $L$ away from a larger neighbourhood of the
boundary.
\end{itemize}

Moreover, $L''$ is independent, up to Hamiltonian isotopy, of the choices which will be made in its construction.
\end{lem}
\begin{proof}
Consider a point $p$ in $\partial L'$, and let $\vec{n}$ be a tangent vector in $TL'_p$ which points towards the interior.  Since $L' $ and $\partial V_{L}$ are transverse, $\vec{n}$ projects to a vector in $TL_p$ which is transverse to the boundary.  Since there is a neighbourhood of $\partial L'$ in $L'$ which is contained in $V_L$, the image of $\vec{n}$ must point towards the interior of $L$. Therefore, the
restriction of the projection $V_L \to L$ to $L'$ is a submersion in a
neighbourhood of the boundary.  In particular, passing to the
cotangent bundle of $L$, we can identify $L'$ locally as the graph of
a closed $1$--form which vanishes on $\partial L$.  Since the inclusion
of $\partial L$ in a neighbourhood induces an isomorphism on
cohomology this $1$--form is exact, so we may write it as the
differential of a function $H'$.

By choosing a cutoff function with appropriately bounded derivatives,
we can construct a function $H'' \co L \to \bR$ with support in a
neighbourhood of $\partial L$ which agrees with $H'$ in a smaller
neighbourhood of the boundary, and such that the graph of $dH''$ lies
in a neighbourhood of $L$ which is identified with $V_L$. Linear
interpolation yields an isotopy between any two choices for $H''$.
Pulling back the graph of the exact $1$--form $dH''$ to $V_L$ yields
the desired Lagrangian $L''$.
\end{proof}

In the situation considered in this paper, $N$ will be a Stein
manifold and all the Lagrangians we will consider will have their
boundary lying on $M$, the zero level set of a holomorphic map $f \co N
\to \bC$ with $0$ as a regular value.  We can equip this Stein
manifold with the structure of an exact symplectic manifold by
choosing an embedding into $\bC^{R}$ and restricting the usual
symplectic form on $\bC^{R}$ and its primitive.  We will denote the
symplectic form by $\omega$, and the primitive by $\theta$.  As usual,
the complex and symplectic structures induce a metric $g$. Given a
Lagrangian $L$, the existence of the complex structure determines
$\partial V_L$ to first order since we can let the line bundle $E$ be
spanned by $J \vec{n}$, with $\vec{n}$ the normal vector of the
inclusion $\partial L \subset L$.

Note that $\theta$ restricts to a closed $1$--form on every Lagrangian
submanifold.  The following definition is standard.
\begin{defin}
A Lagrangian submanifold $L$ is {\em exact} if there is a function $h$
on $L$ such that $dh = \theta|_{L}$.
\end{defin}

In essence, this exactness conditions provides a priori bounds for the
energy of pseudo-holomorphic discs.  Since our Lagrangian submanifolds
may have boundary, we will require an additional condition as a mean
to guarantee some compactness results for such discs (see \fullref{GC}).  Recall that the symplectic orthogonal complement to the tangent space of the fibre defines a distribution on
$N$.  Every tangent vector in $\bC$ has a unique lift to this orthogonal complement, so we may associate a connection to the map $f$.  Given a Lagrangian in the symplectic hypersurface $f^{-1}(0)$ and a curve $\gamma$ in $\bC$ with an endpoint at the origin, parallel transport with respect to this connection determines a unique Lagrangian submanifold of $N$ whose image under $f$ is exactly $\gamma$.

\begin{defin}
A compact oriented exact Lagrangian submanifold $(L, \partial L)$ of
$N$ which has boundary on $M$ is {\em admissible} if there exists a
curve $\gamma$ in $\bC$ such that $\gamma(0)=0$ and $L$ agrees with
the parallel transport of $\partial L$ along $\gamma$ in some
neighbourhood of the origin.
\end{defin}

 Note that the condition of admissibility is of course vacuous if
 $\partial L = \emptyset $, which is a possibility that we do not
 exclude.  As far as the author knows, the idea of studying such
 Lagrangians is due to Kontsevich \acite{kont-ENS}*{page 30}.

As our goal is to define and compute Floer homology for admissible
Lagrangians, we will have to understand the possible behaviours at the
boundary. If $L$ is admissible, then there exists a non-zero vector
$\gamma'(0)=v \in \bR^{2}$ such that the image of a transverse vector
to $T \partial L$ in $TL|_{\partial L}$ lies in the ray $\bR^+ v$.
\begin{defin}
A pair of admissible Lagrangians $(L_1, L_2)$, whose tangent spaces at
the boundary project to vectors $(v_1, v_2)$, is {\em positively
oriented} if the angle from $v_1$ to $v_2$ lies in the interval
$(0, \frac{\pi}{2})$.  If this angle is between
$-\frac{\pi}{2}$ and $0$, we say that the pair is {\em negatively oriented}.
\end{defin}
\begin{figure}
\begin{center}
$\begin{array}{c@{\hspace{.5in}}c} 
\includegraphics[width=2in]{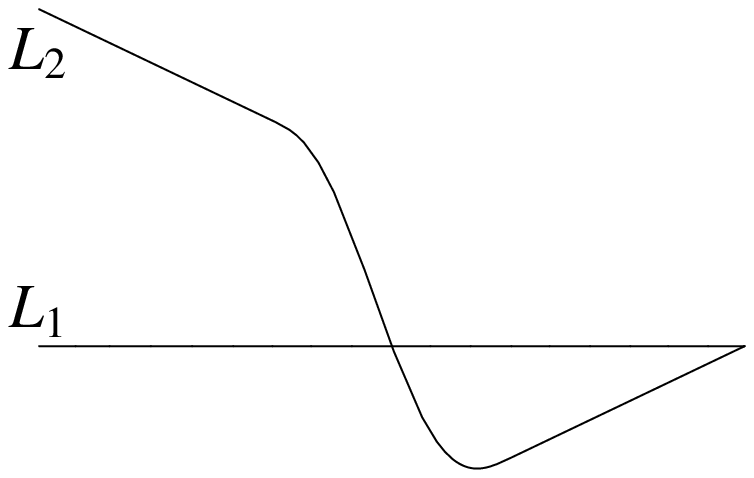} 
& 
\includegraphics[width=2in]{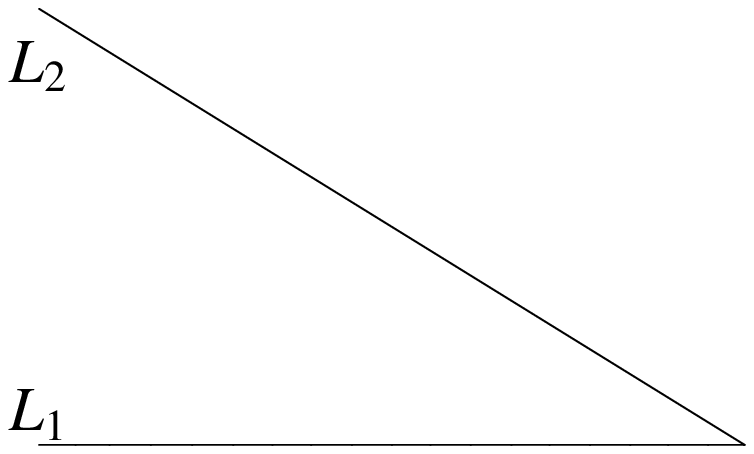}
\end{array}$
\end{center}
\caption{The pair $(L_1, L_2)$ on the left is a positively oriented
pair. The one on the right picture is negatively oriented.}
\label{pos-neg}
\end{figure}

Assume the pair $(L_1, L_2)$ is negatively oriented.  We may choose a
vector $v'_2$ such that $(v_1,v'_2)$ is a positively oriented pair and
the angle between $v_2$ and $v'_2$ is less than $\frac{\pi}{2}$.  The latter
condition is enough to guarantee that the parallel transport of
$\partial L_2$ along a the straight half ray $\gamma'_2$ with tangent
vector $v'_2$ at the origin lies in a neighbourhood of $L_2$ in which
we can apply \fullref{C1-close} to obtain a Lagrangian submanifold
$L_2''$ which interpolates between $L_2$ and the parallel transport of
its boundary along $\gamma'_2$.  This process takes the configuration
on the right in \fullref{pos-neg} to the one on the left. We
conclude the following: 

\begin{lem} \label{change-neg-to-pos}
If the pair $(L_1, L_2)$ is negatively oriented and $H^1(L_2, \partial
L_2) = 0$, then $L_2$ is Hamiltonian isotopic to an admissible
Lagrangian $L''_2$ such that the pair $(L_1, L''_2)$ is positively
oriented.  Furthermore, $L''_2$ is independent of the choices that are made in the construction up to a Hamiltonian isotopy which preserves the admissibility and positivity of the pair $(L_1, L''_2)$.
\end{lem}
\begin{proof}
The only parts which we have not checked are the exactness and
uniqueness of $L''_2$.  But the restriction of $\theta$ to $L_2$, and
therefore to $\partial L_2$ is exact. Since $L_2$ and $L_2''$ are
homeomorphic, the vanishing of $H^{1}(L_2, \partial L_2)$ guarantees
that $\theta|_{L''_2}$ is also exact.

We already proved uniqueness in \fullref{C1-close} by using a linear
isotopy.  But such an isotopy does not necessarily preserve admissibility.  However, any two candidates $L''_{2,0}$ and $L''_{2,1}$  can be made the endpoints of a $1$--parameter family of admissible Lagrangians $L''_{2,t}$ such that $(L_{1}, L''_{2,t})$ is a positively oriented pair for every time $t$.  Since the submanifolds $L''_{2,t}$ agree away from a tubular neighbourhood of the boundary, this $1$--parameter family is generated by a time dependent $1$--form which vanishes away from a neighbourhood of the boundary.  On the other hand, it also vanishes at the boundary since all Lagrangians $L''_{2,t}$ have the same boundary.  Therefore the $1$--form is exact, and the $1$--parameter family $L''_{2,t}$  is generated by a time-dependent Hamiltonian function.
\end{proof}

\subsection{Compactness for pseudo-holomorphic discs}

We will define a Floer theory
for admissible Lagrangians in which the possible boundary intersection
points are ignored.  To this effect, let $L_1$ and $L_2$ be admissible
Lagrangians which intersect transversely away from $M$, and choose $p$
and $q$, a pair of transverse intersection points between these two
Lagrangians.  In analogy with Morse theory, where we count gradient
trajectories connecting critical points, we will count
pseudo-holomorphic maps connecting $p$ and $q$.

Formally, we equip the strip
\[S= \bR \times [0,1] = \{ t , s | - \infty < t < + \infty \, , \, 0 \leq s \leq 1\}\]
with the usual complex structure
\[ J_S \frac{\partial}{\partial t} = \frac{\partial}{\partial s} \, \, \, \, \, \, \, \,  J_S \frac{\partial}{\partial s} = - \frac{\partial}{\partial t} .\]
Just as Morse theory relies on choosing a function whose gradient
flow  with respect to an underlying metric satisfies appropriate transversality conditions, we will have to perturb the complex structure to an almost complex structure in order to achieve transversality, (these almost complex
structures are called regular, Floer--Hofer--Salamon \cite{FHS}).  Perturb the complex
structure on $N$ to such a regular almost complex structure $J$ while maintaining the condition that the map $f \co N \to \bC $ remains holomorphic in a neighbourhood of $M$.  In the closed case, \acite{FHS}*{Remark 5.2} proves that regularity may be achieved by perturbing the almost complex structure on any open subset where every pseudo-holomorphic curve passes.  The case of curves with Lagrangian boundary conditions is entirely analogous, and we can therefore choose a perturbation which is supported in a small neighbourhood of $p$ and $q$, and hence occurs far away from $M$.

Even though the perturbed almost complex structure on $N$ is not
integrable, we may still consider holomorphic maps
\[ u\co S \to N, \]
such that $u(t,0) \in L_1$, $u(t,1) \in L_2$.  We say that such a map
has finite energy if the integral
\[ \int_{S} u^{*}(\omega) \]
is finite.

\begin{figure}
\begin{center}
\labellist\small\hair1pt
\pinlabel $L_1$ [r] at 98 158
\pinlabel $L_2$ [tr] at 136 20
\pinlabel $B_{\epsilon}$ [l] at 360 157
\pinlabel $B_{\epsilon'}$ [l] at 409 157
\endlabellist
\includegraphics[width=2in]{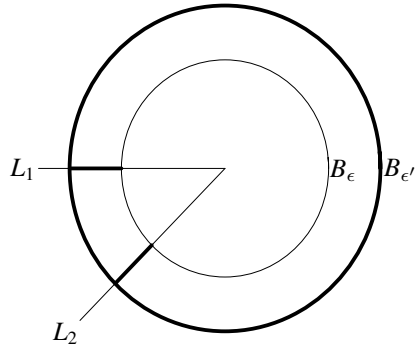} 
\end{center}
\caption{The image of $\partial S'$ is contained in the thick lines.}
\label{comp-fig}
\end{figure}

The classical theory studies $\cM(p,q)$, the moduli space of finite
energy maps such that for every $s \in [0,1]$,
\begin{align*}
\lim_{t \to + \infty} u(t,s) & = p \\ \lim_{t \to - \infty} u(t,s) &
= q.
\end{align*}
We will need to add an additional restriction.
\begin{defin}
Let $L_1$ and $L_2$ be admissible Lagrangians with transverse
intersection points $p$ and $q$ which occur away from the boundary.
An {\em admissible} strip $u$ is an elements of $\cM(p,q)$ whose image
does not intersect $M$.
\end{defin}  
 Note that $\cM(p,q)$ admits a free $\bR$ action corresponding to
translation in the $t$ direction, and the regularity of the almost complex structure implies that the quotient is a smooth manifold.
However, unlike the case for closed Lagrangians, $\cM(p,q)/\bR  $ does not admit a compactification to a manifold with boundary because of those pseudo-holomorphic strips whose boundary intersects the boundaries of $L_1$ or $L_2$.  However, for admissible strips, we have the following:
\begin{lem} \label{GC}
Let $L_1$ and $L_2$ be admissible Lagrangians, and let $p$ and $q$ be
transverse intersection points.  If $u_\tau$ is a $1$--parameter family
of holomorphic strips in $\cM(p,q)$ such that $u_0$ is admissible,
then there exists an $\epsilon >0$ such that the image of $f \circ
u_{\tau}$ does not intersect the closed $\epsilon$ neighbourhood of
the origin for all $0 \leq \tau \leq 1$.  In particular $u_1$ is
admissible.
\end{lem}
\begin{proof}
The reader may find \fullref{comp-fig} useful in what follows.

Let $\epsilon$ be such that $f(L_1)$ and $f(L_2)$ do not intersect in
the punctured $2 \epsilon$ neighbourhood of the origin.  Assume that
$\tau$ is the smallest time at which the lemma does not hold.
Consider $(u_{\tau} \circ f)^{-1} (B_{\epsilon'}(0)) = S'$ for some
$\epsilon'$ between $\epsilon$ and $2 \epsilon$. Since the critical
points of holomorphic maps are isolated, we may choose $\epsilon'$
such that $S'$ is a submanifold of $S$ with boundary.  Notice that the
boundary of $S'$ must be mapped to the union of the set $f(L_1) \cup
f(L_2)$ with the circle of radius $\epsilon'$.

We can compute the degree of $u_{\tau} \circ f$ by choosing a generic
point $p$ in the image of $S'$, and counting the number of preimages.
However, one may find a path from $p$ to the interior of
$B_{\epsilon}(0)$ which does not pass through the image of $\partial
S'$ or through any critical point.  Since the number of preimages can
only change at a critical point or on the boundary, and the image of
$S'$ does not intersect the interior of the ball of radius $\epsilon$,
we conclude that $p$ has no preimages.  But this proves that $u(S)$
does not intersect the ball of radius $\epsilon$.
\end{proof}

Note that the inclusion,
\[ \cM^{\adm}(p,q) \subset \cM(p,q) \]
is clearly open. The above lemma shows that it is also closed, so
$\cM^{\adm}(p,q)$ consists of components of the moduli space of all
finite energy strips connecting $p$ and $q$.   Moreover, the lemma also guarantees that the images of all admissible pseudo-holomorphic strips lie in a
compact subset of $N$ that does not intersect the boundary of the
Lagrangians $L_i$. In particular we can conclude that the closure of $\cM^{\adm}(p,q)/\bR$
in the Gromov compactification $\overline{\cM(p,q) / \bR}$ is a
compact manifold with corners as long as we can prove the existence of a bound on the energy of pseudo-holomorphic strips.  This is where the exactness conditions are used.  Indeed, since $\omega$ is exact, we have
\[ \int_{S} u^{*}(\omega) = \int_{\partial S} \theta \]
by Stokes's theorem.  But $\partial S$ consists of two segments, one on $L_1$ and the other on $L_2$, with endpoints at $p$ and $q$.  Since the restriction of $\theta$ to $L_i$ is exact, the right hand-side is independent of the paths between $p$ and $q$.  This implies that all pseudo-holomorphic strips in $\cM(p,q)$ have the same energy.

\subsection{Floer homology}

Having proved the necessary compactness result, we can now define
relatively graded Floer homology groups over a field of characteristic
$2$ for admissible Lagrangians.  We will follow the construction of
these groups for closed Lagrangians which is due to Floer,
\cite{floer}.  We assume $(L_1 ,L_2)$ is a positively oriented pair of
Lagrangians which intersect transversely away from the boundary, and
define a chain complex
\[ CF_{*}(L_{1}, L_{2}) = \bigoplus_{p \in (L_{1} \cap L_2) - M} \bZ_{2} \cdot [p] \]
with differential
\[ d[p] = \sum_{[q]} | \cM^{\adm}(p,q) / \bR| \cdot [q] \]
where the sum is taken over all points $q$ such that $\cM^{\adm}(p,q)$
is $1$--dimensional, and $|\cM^{\adm}(p,q) /\bR|$ is the cardinality of
the space of unparametrized strips connecting $p$ to $q$.  The
following result is classical.
\begin{lem}
If $L_1$ and $L_2$ are exact, then $d^{2}=0$.
\end{lem}
\begin{proof}
The proof that $d^{2}=0$ relies on interpreting the terms in the
expression for $d^2 [p]$ as the boundaries of $1$--dimensional moduli
spaces of unparametrized strips.  So we must show the boundary of such
moduli spaces consists only of pairs of holomorphic strips, ie,\ that
no bubbling of pseudo-holomorphic discs occurs.  But given a
pseudo-holomorphic disc $u$ with boundary on one of the Lagrangians
$L_i$ we can use Stokes's theorem to compute that
\[ 0 \neq \int_{D} u^{*}(\omega) = \int_{S^1} u^{*}(\theta).\]
This contradicts the assumption that the restriction of $\theta$ to
$L_i$ is exact.
\end{proof}

This construction produces a $\bZ_2$ graded theory which cannot be
lifted, in general, to the usual $\bZ$ grading which we expect in a
homology theory.  However, if $N$ admits a complex volume form (ie,\
if $c_1(N)=0$), then there is a special class of Lagrangians for
which such a theory exists, see Seidel \cite{seidelGL}. Choosing a complex volume
form $\Omega$, we define a phase map
\begin{align*}
L & \to  S^{1} \\ p & \mapsto  \frac{\Omega(e_1 \wedge \cdots \wedge
e_n)}{\left| \Omega( e_1 \wedge \cdots \wedge e_n ) \right|}
\end{align*}
where  $\{ e_{i}
\}_{i=1}^{n}$ is an oriented frame for the tangent space of $L$ at
$p$.

\begin{defin}
A {\em graded Lagrangian} submanifold of $N$ is an oriented Lagrangian
submanifold of $N$ together with a lift of its phase map from $S^{1}$ to $\bR$.
\end{defin}

Note, in particular, that the obstruction to the existence of a lift
lies in $H^{1}(L)$, so that all simply connected Lagrangians are
gradable.  We refer to \cite{seidelGL} for details on the construction
of $\bZ$--graded Floer homology.  We also need to lift these homology
groups to $\bC$--valued invariants in order to compare them to
cohomology groups of line bundles on the mirror.  The need for
appropriate orientations accounts for the additional restrictions in
\begin{defin}
An {\em admissible Lagrangian brane} is an admissible graded
Lagrangian $L$ which is spin, together with a choice of a spin
structure.
\end{defin}
\begin{rem}
The term brane is borrowed from string theory, where the ``Lagrangian
branes'' that we're considering can be thought of as boundary
conditions for open strings in the $A$--model.  Other than the fact
that string theory motivates the homological mirror symmetry
conjecture, physical considerations are completely irrelevant to our
arguments.
\end{rem}

\begin{figure}
\begin{center}
\begin{picture}(0,0)%
\includegraphics{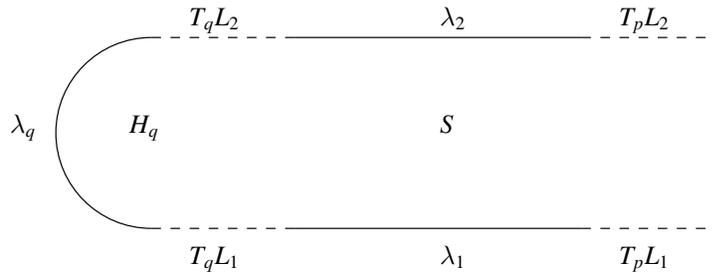}%
\end{picture}%
\setlength{\unitlength}{3947sp}%
\begingroup\makeatletter\ifx\SetFigFont\undefined%
\gdef\SetFigFont#1#2#3#4#5{%
  \reset@font\fontsize{#1}{#2pt}%
  \fontfamily{#3}\fontseries{#4}\fontshape{#5}%
  \selectfont}%
\fi\endgroup%
\begin{picture}(4778,1714)(2101,-3044)
\put(3226,-1486){\makebox(0,0)[lb]{\smash{{\SetFigFont{12}{14.4}{\familydefault}{\mddefault}{\updefault}{\color[rgb]{0,0,0}$T_{q}L_2$}%
}}}}
\put(3226,-2986){\makebox(0,0)[lb]{\smash{{\SetFigFont{12}{14.4}{\familydefault}{\mddefault}{\updefault}{\color[rgb]{0,0,0}$T_{q}L_1$}%
}}}}
\put(4801,-1486){\makebox(0,0)[lb]{\smash{{\SetFigFont{12}{14.4}{\familydefault}{\mddefault}{\updefault}{\color[rgb]{0,0,0}$\lambda_2$}%
}}}}
\put(4801,-2986){\makebox(0,0)[lb]{\smash{{\SetFigFont{12}{14.4}{\familydefault}{\mddefault}{\updefault}{\color[rgb]{0,0,0}$\lambda_1$}%
}}}}
\put(5926,-1486){\makebox(0,0)[lb]{\smash{{\SetFigFont{12}{14.4}{\familydefault}{\mddefault}{\updefault}{\color[rgb]{0,0,0}$T_{p}L_2$}%
}}}}
\put(5926,-2986){\makebox(0,0)[lb]{\smash{{\SetFigFont{12}{14.4}{\familydefault}{\mddefault}{\updefault}{\color[rgb]{0,0,0}$T_{p}L_1$}%
}}}}
\put(2101,-2161){\makebox(0,0)[lb]{\smash{{\SetFigFont{12}{14.4}{\familydefault}{\mddefault}{\updefault}{\color[rgb]{0,0,0}$\lambda_q$}%
}}}}
\put(2851,-2161){\makebox(0,0)[lb]{\smash{{\SetFigFont{12}{14.4}{\familydefault}{\mddefault}{\updefault}{\color[rgb]{0,0,0}$H_q$}%
}}}}
\put(4801,-2161){\makebox(0,0)[lb]{\smash{{\SetFigFont{12}{14.4}{\familydefault}{\mddefault}{\updefault}{\color[rgb]{0,0,0}$S$}%
}}}}
\end{picture}%
\end{center}
\caption{Gluing $H_q$ to $S$}
\label{glueing}
\end{figure}

The relevance of spin structures to the orientation of moduli spaces
of discs was observed by de Silva in \cite{deSilva}  and by Fukaya, Oh, Ohta, and Ono in \cite{FOOO3}.  We will give a short description of the signed differential.  The reader should keep \fullref{glueing} in mind during the next few paragraphs.

Consider an intersection point $p$ between $L_1$ and $L_2$.  The
tangent spaces $T_p L_1$ and $T_p L_2$ are linear Lagrangians in a
symplectic vector space $V$.  The grading determines a unique path up
to homotopy $\lambda_p$ from $T_p L_1$ to $T_p L_2$. We may therefore
define a Cauchy--Riemann operator $\bar{\partial}_{H_p}$ on a copy of
the upper half-plane $H_{p}$ with Lagrangian boundary conditions, given
by the path $\lambda_p$, that converge to $L_1$ along the negative real
axis and to $L_2$ along the positive real axis.

More precisely, choose a map $h\co \bR \to [0,1]$ which takes a
neighbourhood of $-\infty$ to $0$ and of $\infty$ to $1$, and consider
maps from $H_p$ to the vector space $V$ that have values in the
Lagrangian $\lambda_p (h(t))$ subspace at a point $t$ of the boundary.  We
denote this space by $\cC_{\lambda_p}^{\infty}(H_{p},V)$.  We can now
define the operator
\[ \bar{\partial}_{H_p}\co \cC_{\lambda_p}^{\infty}(H_{p},V) \to \cC^{\infty}(H_{p}, V \otimes \Omega^{0,1}(H_p)) \]
to be the usual $\bar{\partial}$ operator.  We can perform the same
construction at any other intersection point $q$ to produce an
operator $ \bar{\partial}_{H_q}$ with boundary conditions $\lambda_q$.

Let $S$ be an element of $\cM^{\adm}(p,q)$.  Since the upper half-plane
is biholomorphic to a $1$--sided strip, we may glue $H_q$ to the
negative end of $S$ to yield a surface $S \# H_q$ which is again
bi-holomorphic to a one-sided strip.  Strictly speaking, this means that we should choose an identification of a neighbourhood of infinity in $H_q$ with $(-\infty, 0] \times [0,1]$.  For any $R >0$ we can remove $ (-\infty, -R)  \times [0,1] $ from $S$ and from $H_q$, then identify the two copies of $ [-R, 0] \times [0,1] $ in $S$ and in $H_q$ using the biholomorphism
\[ (t,s) \mapsto (-R -t, 1 -s) .\]
Further, for $R$ sufficiently large, since the asymptotic
boundary conditions agree, we can glue the Cauchy--Riemann operators
$\bar{\partial}_{S}$ and $\bar{\partial}_{H_q}$ to get a
Cauchy--Riemann operator $\bar{\partial}_{S} \# \bar{\partial}_{H_q}
$on the strip $S \# H_q$.  We do not keep track of $R$ in the notation because the construction is essentially independent of $R$ if it is sufficiently large.

To make sense of the boundary conditions, we trivialize the restriction
of $TN$ to $S$, so that the boundary conditions $u(t,0) \in L_1$ and
$u(t,1) \in L_2$ yield two paths
\begin{align*}
\lambda_{1}\co \bR & \to  V \textrm{ such that } \lambda_{1}(-\infty)
= T_{q}L_1 \textrm{ and } \lambda_{1}(+ \infty) = T_{p}L_1 \\
\lambda_{2}\co \bR & \to  V \textrm{ such that } \lambda_{2}(-\infty)
= T_{q}L_2 \textrm{ and } \lambda_{2}(+ \infty) = T_{p}L_2.
\end{align*}
The boundary condition $\lambda_{S \# H_q}$ for $\bar{\partial}_{S} \#
\bar{\partial}_{H_q} $ is just given by concatenating the paths
$-\lambda_1$, $\lambda_q$, and $\lambda_2$, where $-\lambda_1$
corresponds to traversing the path ``backwards''.  In particular, we
have the asymptotic conditions
\[ \lambda_{S \# H_q} ( - \infty) = T_{p}L_1 \textrm{ and }  \lambda_{S \# H_q} (+\infty) = T_{p}L_2 \]
which are the same conditions as those of $\lambda_p$.  Since the
space of graded Lagrangians is simply connected, the paths $\lambda_p$
and $\lambda_{S \# H_q}$ are homotopic.  Moreover, the choice of spin
structures determines, up to homotopy, a unique such homotopy between  $\lambda_p$ and $\lambda_{S \# H_q}$.  If we
consider the extensions of our Cauchy--Riemann operators to appropriate
Hilbert space completions of spaces of smooth functions, the above
path will therefore yield an isomorphism
\[   \det(\bar{\partial}_{S \# H_q}) \cong \det(\bar{\partial}_{H_p}). \]
On the other hand, the gluing theorem yields the following isomorphism
of determinant bundles
\[ \det(\bar{\partial}_{S \# H_q}) \cong \det(\bar{\partial}_{S}) \otimes \det(\bar{\partial}_{H_q}).\]
If $S$ is a strip with $1$--dimensional parametrized moduli space, then
$\det(\bar{\partial}_{S})$ is canonically trivialized by the
translation operator $\frac{\partial}{\partial t}$.  We therefore
obtain an isomorphism
\[ \det(\bar{\partial}_{H_p}) \cong \det(\bar{\partial}_{H_q}).\]
The orientations of $L_1$ and $L_2$ determine orientations of these two vector spaces.  The contribution of $S$ to $d[p]$ will be positive if the above
isomorphism preserves orientations, and negative otherwise.

\subsection{Stability of Floer homology under perturbations}
While we are interested in the symplectic topology of a complex
hypersurface $M \subset \ctorus{n}$, our constructions will rely on
deforming $M$ to symplectic submanifolds of $\ctorus{n}$ that are not
necessarily complex.  We will therefore have to generalize the above
discussion to the non-integrable situation. There is a more general
setting of exact symplectic manifolds with properly embedded
hypersurfaces to which this discussion can be extended, but our aims
here are more modest.  We remark that our construction of Floer
homology for admissible Lagrangians is still valid if we perturb both
the complex structure and the fibration near the $0$--level set so long
as the following conditions are preserved:
\begin{itemize}
\item The perturbed almost complex structure is compatible with the
symplectic form and $J$--convex at infinity.
\item The map to $\bC$ near the $0$--level set remains holomorphic.
\end{itemize}
The first condition is familiar from the study of the symplectic
topology of Stein manifolds Eliashberg--Gromov
\cite{eliashberg-gromov}.  The second
condition is necessary for the validity of \fullref{GC}.

Note that admissibility is not stable under exact Hamiltonian
perturbations (even those that preserve $M$).  However, given two
positively oriented admissible Lagrangians $L_1$ and $L_2$ one may
still consider exact Hamiltonian deformations which preserve the
admissibility of $L_1$. As long as the tangent vectors $v_1$ and $v_2$
are appropriately oriented, \fullref{GC} will apply for some
$\epsilon$.  In particular, the usual proofs of invariance of Floer
homology are valid for those Hamiltonian isotopies that preserve the
admissibility of $L_1$, and the positivity of the pair $(L_1, L_2)$.
As in the classical situation, this allows us to compute Floer
homology for a pair of admissible Lagrangian branes that do not
intersect transversely away from $M$ by choosing an appropriate
Hamiltonian deformation of one of them.

In particular, if $(L_{1}, L_{2})$ is a pair of admissible Lagrangian
branes which is negatively oriented, we can use \fullref{change-neg-to-pos} in order to unambiguously define
\[ CF_{*}(L_{1}, L_{2}) \equiv CF_{*}(L_{1}, L''_{2})\]
where the pair $(L_1, L''_2)$ is positively oriented.

We also observe that our conventions for ``positivity'' are designed
to guarantee that Floer's old result that expresses Floer homology of
nearby Lagrangians in terms of ordinary homology extends to this
setting,
\[ HF_{*}(L, L) \cong H_{*}(L, \partial L) .\]
We now consider what happens when we vary the holomorphic map $f$ in a
family.
\begin{lem}\label{invariance}
Let $f_{t}$ be a $1$--parameter family of symplectic fibrations which
are holomorphic near the origin, and assume that
$(f_t)^{-1}(0)=M$. There is a bijection between Hamiltonian isotopy
classes of admissible Lagrangians for $f_0$ and $f_1$. For Lagrangians
satisfying $H^1(L,\partial L) =0$, this bijection respects the Floer
homology groups.
\end{lem}
\begin{proof}
Assume that $L$ is an admissible Lagrangian with respect to $f_0$.
Since $\partial L$ is compact, we can uniformly bound the derivatives
$\frac{\partial f}{\partial t}$ in a neighbourhood of $\partial L$.
We can therefore choose an $\epsilon$ such that the parallel transport
of $\partial L$ with respect to $f_{t}$ lies in an arbitrarily small
neighbourhood of the parallel transport of $\partial L$ with respect
to $f_{t+\epsilon}$.  If this neighbourhood is small enough, we can
apply \fullref{C1-close}.  In particular, subdividing the interval
$[0,1]$ into sufficiently many subintervals, we obtain a Hamiltonian
isotopy between $L$ and a Lagrangian which is admissible with respect
to $f_1$.

Let $L_1$ and $L_2$ be two admissible Lagrangians with respect to
$f_0$.  We can choose the Hamiltonian isotopies which we used in the
previous paragraph to be supported in an arbitrarily small
neighbourhood of $M$, and such that no new intersection points are created.  In particular, these isotopies will equal the identity near the images of the elements of $\cM(p,q)$ for all pairs of intersection $p$ and $q$ between $L_1$ and $L_2$.  This yields the desired invariance of Floer homology groups.
\end{proof}

\subsection{Cohomology and cup product}
With Poincar\'e duality in mind, we can now define Floer cohomology
by simply re-indexing the Floer complex as
\[ CF^{*}(L_1,L_2) \equiv CF_{n-*}(L_1,L_2)\]
which on cohomology yields
\[ HF^{*}(L_1,L_2) \equiv HF_{n-*}(L_1,L_2).\]
\begin{rem}
Note that Floer homology and cohomology are Poincar\'e dual with these conventions, but there is no degree preserving duality between them.  Whereas the classical analogue of our Floer homology groups is the homology of $L$ relative its boundary, the classical analogue of our Floer cohomology group is the ordinary cohomology of $L$.  These classical groups are indeed Poincar\'e dual.  One can resolve this unfortunate state of affairs by introducing Floer homology in two flavours, with positive and negative orientations at the boundary.  In this language, our Floer cohomology between $L_1$ and $L_2$ would indeed be the (ordinary) dual of Floer homology between $L_2$ and $L_1$ with the opposite convention to the one we have chosen.  Since the main goal of this paper is to perform a computation, we will not discuss these issues further, and simply use the above definition.
\end{rem}
As in the case of closed Lagrangians the Floer co-chain complex admits
a cup product
\[ CF^{*}(L_1, L_2) \otimes CF^{*}(L_2,L_3) \to  CF^{*}(L_1,L_3) .\]
In order to describe the cup product, we need some preliminary
definitions.  There are essentially no differences between the case we
are studying, and that of closed Lagrangians which is explained in
Fukaya and Oh \cite{FO}.
\begin{defin}
A Riemann surface with {\em strip-like ends} is an open Riemann
surface with a choice of biholomorphisms between its ends (ie,\
complements of sufficiently large compact subsets) and the strip
$[0,1] \times \bR^{+}$.
\end{defin}

\begin{figure}
\begin{center}
\includegraphics[width=2in]{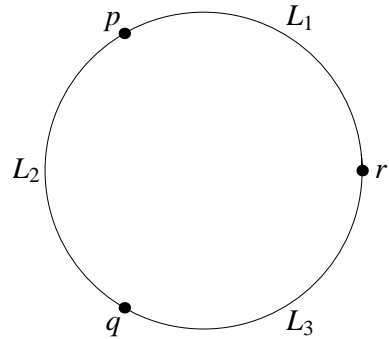}
\end{center}
\caption{A marked disc}
\label{labled_circle}
\end{figure}

We let $D$ be the unit disc in $\bC$, and $\xi = e^{\frac{2 \pi
i}{3}}$.  Note that $T=D - \{\xi, \xi^{2}, 1\}$ admits the structure
of a Riemann surface with strip-like ends. Assuming all Lagrangians
intersect transversely and the pairs $(L_1,L_2)$, $(L_2,L_3)$ and
$(L_1,L_3)$ are positively oriented, we consider finite energy
pseudo-holomorphic maps $u\co T \to N$ which satisfy the following
conditions (see \fullref{labled_circle}):
\begin{itemize}

\item $u$ maps the arcs $(1, \xi)$ to $L_1$, $(\xi, \xi^{2})$ to
$L_2$, and $(\xi^{2}, 1)$ to $L_3$.

\item Along the strip-like ends associated to the punctures $(\xi,
\xi^{2}, 1)$, the image of $u$ converges uniformly to interior intersection
points $(p,q,r)$ among the Lagrangians $L_j$.

\item The image of $T$ under $u$ does not intersect $M$.

\end{itemize}

We denote the above moduli space by $\cM^{\adm}(p,q,r)$.  We can now
define the cup product over $\bZ_{2}$ by the formula
\[ [p] \otimes [q] \mapsto \sum_{r}  |\cM^{\adm}(p,q,r)| \cdot [r] ,\]
where the sum is taken over $0$--dimensional moduli spaces.  Just as we
proved that $d^{2}=0$, the usual proof that this cup product descends
to an associative product on cohomology applies in our
situation.

To obtain a product in Floer cohomology over $\bC$, we follow the same
strategy as for obtaining signs in the differential.  Given a
holomorphic triangle $T$ with strip like ends associated to three
points $p$, $q$ and $r$, we attach copies of the upper half-plane at
the incoming ends to form a $1$--sided strip $T \# H_p \# H_q$.  As
before, we obtain an operator $\bar{\partial}_{T \# H_p \# H_q }$.
The gluing theorem yields a canonical isomorphism
\[ \det(\bar{\partial}_{T \# H_p \# H_q  }) \cong \det(\bar{\partial}_{T}) \otimes \det(\bar{\partial}_{H_p}) \otimes \det(\bar{\partial}_{H_q}) .\]
There is a unique path (up to homotopy) between the
boundary conditions for this operator and those for the operator
$\bar{\partial}_{H_r}$, yielding a canonical isomorphism
\[ \det(\bar{\partial}_{T \# H_p \# H_q  }) \cong \det(\bar{\partial}_{H_r}) .\]
Since $u$ only contributes to the cup product when its associated
$\bar{\partial}$ operator is invertible, we have a canonical
trivialization of $\det(\bar{\partial}_{T})$. So we obtain an
isomorphism
\[  \det(\bar{\partial}_{H_p}) \otimes \det(\bar{\partial}_{H_q}) \cong \det(\bar{\partial}_{H_r}) ,\]
whose compatibility with orientations determines the signed
contribution of $T$ to the image of $[p] \otimes [q]$ under the cup
product.

\section{Tropical geometry} \label{TG}
\subsection{Background}
Let $f$ be a Laurent polynomial in $n$ variables over $\bC$.  In
\cite{GKZ}, Gel'fand, Kapranov, and Zelevinsky introduced the
amoeba of $f$ as the projection of $f^{-1}(0)$ to $\bR^{n}$ under the
logarithm map
\[ \Log(z_1, \ldots, z_{n}) = (\log|z_1|, \ldots , \log|z_{n}| ) .\]
One may also define amoebas for varieties over fields other than
$\bC$.  Amongst other results on varieties over non-Archimedean
fields, Kapranov gave a description of their amoebas in term of
polyhedral complexes in $\bR^{n}$ \cite{EKL}.  Mikhalkin then used
this description in order to obtain new results about the topology of
complex hypersurfaces \cite{mikhalkin}.  We will follow his point of
view with small modifications.

In coordinates, we will write \[ f = \sum_{\alpha \in \bZ^{n}} c_{\alpha} z^{\alpha}\]
with $c_{\alpha} \in \bC$ and $z^{\alpha} = z_{1}^{\alpha_1} z_{2}^{\alpha_2} \cdots z_{n}^{\alpha_n}$.  Since $f$ is a polynomial, only finitely many monomials have non-zero coefficients, and we let $A$ be the set of vectors which label these monomials.  Let $P$ be the convex hull of $A$ thought of as a subset of $\bR^{n}$ (Newton polytope of $f$).

\begin{defin}
The {\em convex support} of a function $\nu\co A \to \bR$ is the largest
convex function $\hat{\nu}\co P \to \bR$ such that $\hat{\nu}(\alpha)
\leq \nu(\alpha)$.
\end{defin}

The function $\hat{\nu}$ is in fact piece-wise linear on $P$, so we
may decompose $P = \bigcup P_{\nu_i}$, where each $\nu_i$ is a linear
function, and $P_{\nu_i}$ is the domain where $\nu_i$ and $\hat{\nu}$
agree.  We say that the decomposition $P = \bigcup P_{\nu_i}$ is the
{\em coherent subdivision} of $P$ induced by $\nu$.  The following
result is well known.
\begin{lem}
Each polytope of the coherent subdivision induced by a function $\nu\co
A \to \bR$ is a lattice polytope whose vertices lie in $A$. \noproof
\end{lem}

We will be particularly interested in subdivisions which satisfy the
following additional condition (See \fullref{th}).
\begin{defin}
A subdivision of $P$ is {\em maximal} if each polytope of the
subdivision is equivalent under the action of $ASL(n,\bZ)$ to the
standard $n$--dimensional simplex.
\end{defin}

The function $\nu$ also determines a piecewise linear function
$L_{\nu}\co \bR^{n} \to \bR$, its Legendre transform, which is defined
by,
\begin{equation} \label{legendre} L_{\nu}(u) = \max_{\alpha \in A} \left( \langle \alpha, u \rangle - \nu(\alpha) \right) .\end{equation}
Since it is defined as the maximum of finitely many linear functions,
$L_{\nu}$ is smooth away from a subset of zero measure.

\begin{defin}
Given any function $\nu\co A \to \bR$, the locus of non-smoothness of $L_{\nu}$ is a {\em tropical
hypersurface} or {\em tropical amoeba} of $f$.
\end{defin}

Note that this tropical hypersurface of $f$ is the set where two or
more linear functions agree.  It will be sufficient for our purposes
to think of the tropical amoeba $\Pi$ as an $n-1$ dimensional
polyhedral complex whose $k$ skeleton corresponds to the set where
$n-k+1$ functions agree. In particular, each $k$--face is dual to a
unique $(n-k)$--polytope of the polyhedral subdivision of $A$ induced
by $\nu$.  We will denote the dual of $\tau$ by $\check{\tau}$.  Given a face $\sigma$ of $\Pi$, we will also use $\check{\sigma}$ for its dual polytope in $P$. Our conventions are that all polytopes and faces are closed.

When $k=0$, we obtain a component of $\bR^{n} - \Pi$ rather than a
face of $\Pi$.  Concretely, every component of $\bR^{n} - \Pi$ may be
labeled by the unique $\alpha \in A$ which, on the given component,
achieves the maximum in the equation defining $L_{\nu}$.  We can
therefore write
\[ \bR^{n} - \Pi = \bigsqcup_{\alpha \in A} C_{\alpha} - \partial C_{\alpha} .\]
In general, some of these components may be empty.

Going back to the complex numbers, we may use the function $\nu$ to
define the patchworking polynomials
\[ f_{t}(z) = \sum_{\alpha \in A} c_{\alpha} t^{-\nu(\alpha)} z^{\alpha} .\]
For every $t$, we consider $\cA_{t}$ the $\Log$--amoeba of $f_t$. One
expression of the connection between tropical and complex amoebas is
given by the following theorem which is due to Mikhalkin
\acite{mikhalkin}*{Theorem 5} and Rullg{\aa}rd
\acite{rullgaard}*{Theorem 9}.

\begin{thm} \label{GH}
In the Gromov--Hausdorff topology, the sets $\cA_t / \log(t)$ converge
to $\Pi$ as $t$ goes to infinity. \noproof
\end{thm}

We will be studying the hypersurface $f^{-1}(0) \subset \ctorus{n}$
from the symplectic point of view.  In order to be precise, we must
choose a symplectic structure on $\ctorus{n}$.  We will use
\[ \omega = \sum_{j=1}^{n} \frac{dz_{j} \wedge d\bar{z}_j}{2 i \left|z_{j}\right|^2}, \] 
where $\{z_j\}_{j=1}^{n}$ are the standard coordinates on
$\ctorus{n}$.  With respect to the usual $\ctorus{n}$ action, this is
an invariant K\"ahler form on $\ctorus{n}$.  The vectors $\{
\frac{\partial}{\partial z_{j}} \}_{j=1}^{n}$ and their complex
conjugates form an orthogonal basis for the K\"ahler metric. In
particular,
\[ \left| \frac{\partial}{\partial z_{j}} \right| = \frac{1}{\left|z_{j}\right|} \]
and hence
\[ \left|dz^{j}\right| = \left|z_{j}\right| .\]
It will be convenient to have a different description of $\ctorus{n}$
as a symplectic manifold.  Consider the cotangent bundle of $\bR^n$
with standard coordinates $\{ u_{j} \}_{j=1}^{n}$ on $\bR^{n}$ and $\{
\theta_{j}=du_j \}_{j=1}^{n}$ on the fibre of the cotangent bundle.
The coordinates $\{ u_{j} \}_{j=1}^{n}$ induce an affine structure on
the base, which determines a lattice in the fibre.  We will choose
this lattice to be spanned by $\{ 2 \pi \theta_{j} \}_{j=1}^{n}$.
\begin{lem}
 The quotient of $T^{*} \bR^{n}$ by the lattice $2 \pi \bZ^{n}$ in each cotangent fibre is
 symplectomorphic to $\ctorus{n}$ with the identification given by the
 exponential map
$$ ( u_{j}, \theta_{j} ) \to e^{u_{j} + i \theta_j} .\eqno{\qed}$$
\end{lem}
Note that the inverse of the above exponential map is the logarithmic
map whose first component arises in the definition of amoebas,
\[ z_{j} \to (\log|z_j|, \arg(z_j) ) .\]
We will omit all notation for the above symplectomorphisms, and will
use either coordinate system at our convenience.  Note, however, that
these symplectomorphisms identify the projection of the cotangent
bundle onto its base manifold with the $\Log$ projection of
$\ctorus{n}$ onto $\bR^{n}$, hence identify the cotangent fibres with
the set of points in $\ctorus{n}$ whose components have fixed norms.  We will also be using the standard Euclidean metric on
$\bR^{n}$.  As usual this metric identifies the tangent and cotangent
bundles, and carries the natural complex structure of the tangent
bundle to a complex structure on the cotangent bundle which is
compatible with the symplectic form.

The advantage of the cotangent bundle point of view is that many
constructions can now be performed on the base $\bR^{n}$, and some of
them reduce to linear algebra.  For example, recall that a
diffeomorphism of a manifold induces a symplectomorphism of its
cotangent bundle.  If, in addition, this diffeomorphism is an
isomorphism of the affine structure, it will induce a
symplectomorphism of the associated Lagrangian torus bundle.  Thus,
every affine transformation in $ASL(n,\bZ)$ induces a
symplectomorphism of $\ctorus{n}$.  Assume, for simplicity, that we
have a linear transformation.  If we represent it by a matrix $A$,
then in the standard coordinates of $T^{*}\bR^{n}$ this
symplectomorphism is given by multiplying by $A$ in the coordinates
corresponding to the base, and by $(A^{T})^{-1}$ in the coordinates
corresponding to the fibre.  In general such a transformation does not
preserve the standard Euclidean metric on the base, hence does not
preserve the complex structure on the fibre.

Indeed, given an element of $ASL(n,\bZ)$, there is a different
automorphism of the torus fibration which does preserve the complex
structure.  Assuming again that we have a linear transformation
represented by an integral matrix $A$, this complex automorphism, in
the standard coordinates of $T^{*}\bR^{n}$, is given by using $A$ in
both the base and the fibre directions.  In the standard coordinates
of $\ctorus{n}$ as a complex manifolds, it can be thought of as the
multiplicative change of variables
\[ (z_1, \ldots , z_n) \to (z^{\alpha_1}, \ldots, z^{\alpha_n}) \]
where $\alpha_i$ are the rows of $A$.

\subsection{The tropical model}

Let $X$ be a complete toric variety.  Let $v_{i}$ be the primitive
vertices of the $1$--cones of the fan $\Delta$ defining $X$.  We
consider $A = \{ v_{i} \}_{i=1}^{m} \cup \{0\}$ and its convex hull
$P$.  To each maximal cone $\tau \in \Delta(n)$ we assign the convex
polytope $P(\tau)$ which is the convex hull of the set $\{0 \} \cup
\{ v_{i} | v_{i} \in \tau \}$.

We will assume that every polytope $P(\tau)$ is a minimal simplex;
this is equivalent to the requirement that the set $\{ v_{i} | v_{i} \in
\tau \}$ forms a basis for $\bZ^{n}$. Fans with this property
correspond to smooth toric varieties, \acite{fulton}*{page 34}.  Let $W$
be the Laurent polynomial
\[ W(z) = -1+ \sum_{0 \neq \alpha \in A} z^{\alpha} \] and let $M = W^{-1}(0)$.
\begin{rem}
The choice of coefficients on the monomials which appear in $W$ is only done for convenience.  After passing to the tropical limit, any other choice of (non-zero) coefficients is related to the one we picked by rescaling followed by a symplectomorphism.
\end{rem}

\begin{lem}
Every ample line bundle on $X$ induces a coherent subdivision of $P$
which near $0$ is given by the polytopes $P(\tau)$.  In particular,
near $0$, we have a maximal triangulation.
\end{lem}
\begin{proof}
Every ample line bundle on $X$ is determined up to an integral linear
function by a strictly convex function $\phi$ on $|\Delta|$ which is
integral linear on each cone.  Consider the coherent subdivision of
$P$ induced by $\phi|_{A}$.  To prove the lemma, it will suffice to
prove that $0$ is a vertex of this subdivision, and that each
full-dimensional polytope of the subdivision abutting on the origin is
equal to $P(\tau)$ for some cone $\tau \in \Delta(n)$.

To prove this we observe that on each polytope $P(\tau)$ the convex
support $\psi$ of the function $\phi|_{A}$ agrees with $\phi$ since
the latter is convex.  Since $\phi$ is in fact strictly convex, we
conclude that each $P(\tau)$ is contained in a distinct polytope of
the subdivision.  It remains to show that $P(\tau)$ is equal to a
polytope of the subdivision.  Indeed, if this were not the case, this
polytope would have a vertex not belonging to $\tau$.  But coherent
subdivisions do not create new vertices, so this vertex is in fact the
minimal vertex $v_{\tau'}$ of some $1$--cone $\tau'$.  This contradicts
the strict convexity of $\phi$.
\end{proof}

The choice of an ample line bundle therefore induces a choice of patchworking polynomials $W_{t}$, whose amoebas converge (after rescaling) to a tropical hypersurface $\Pi$ as in the previous section.  Since the origin is a vertex of the subdivision induced by $\phi$, we
conclude that there must be a non-empty component, $Q$, of the
complement of $\Pi$ which is dual to the origin.  

\begin{cor}
$Q$ is the convex polytope which, in toric geometry, describes $X$ as
a toric variety with ample line bundle $\cO(1)$.
\end{cor}
\begin{proof}
Observe that Equation \eqref{legendre} shows that $Q$ is the polytope
consisting of points $y$ which satisfy
\[ \langle v_{i}, y \rangle \leq \phi(v_i)\]
for every primitive vertex $v_i$.  After adjusting for different sign
conventions, this is therefore the polytope whose integral points form
a basis for the space of sections of the ample line bundle $\cO(1)$
\acite{fulton}*{Section 3.4}.
\end{proof}

We do not get a global triangulation of our Newton polytope, but we
certainly have a maximal triangulation near the origin.  Since we will
only be studying $\Pi$ in a neighbourhood of $Q$, we may appeal to the
results of Mikhalkin about maximal tropical degenerations.

\begin{rem} \label{int}
While we could use any ample line bundle, it will be convenient to
have a lattice point in the interior of $Q$.  Note that this can be
achieved by replacing $\cO(1)$ by a sufficiently high power, and
we may assume that the origin is an interior point.
\end{rem}

\subsection{Twisting the tropical zero-section} \label{twist1}
We begin by producing a ``tropical" version of our construction.  In
particular, the boundary of this version will not be smooth, hence it
will not be clear how the results that we obtain would be invariant
under a natural class of Hamiltonian isotopies as discussed in the
previous  section.  Nonetheless, as the complete version requires many
choices that are necessary to produce a meaningful answer, but that
obscure the simple nature of the construction, we will prefer to
discuss the tropical case first.

We now consider the situation where our toric variety $X$ is smooth.
In this situation, we may restate a result of Mikhalkin.
\begin{thm} \acite{mikhalkin}*{Lemma 6.2} \label{th}
Given a maximal triangulation of $P$, there exists a natural choice of
a piecewise smooth symplectic hypersurface $M_{\infty}$ of
$\ctorus{n}$ which projects to the tropical amoeba $\Pi$.  \noproof
\end{thm}
In fact, we only need this theorem as motivation, since we will be using \fullref{family} to prove the precise results.  We will
therefore not give a complete description of $M_{\infty}$ which may be
thought of as a limit of $W_{t}^{-1}(0)$ after rescaling, but we will
use the following results that follow from the proof of \fullref{th} or of \fullref{family}.
\begin{itemize}
\item The preimage of a point on a $k$--face $\sigma$ of $\Pi$ contains
a subtorus of the fibre $\bR^{n}/ 2\pi \bZ^n$ which is parallel to the
tangent space of $\sigma$ (thought of modulo $\bZ^n$ of course).
\item The preimage a point on the interior of an $n-1$ dimensional
facet $\sigma$ is equal to a torus in the fibre which is parallel to the
tangent space of $\sigma$.
\end{itemize}

\begin{rem} \label{Tsk}

We have chosen the coefficients of the monomials that appear in $W$
exactly in such a way as to ensure that the intersection of
$M_{\infty}$ with the zero section contains $\partial Q \subset \Pi$,
the boundary of the component of $\bR^{n} - \Pi$ which corresponds to
the origin.  This guarantees that the preimage of an $n-1$ facet of
$Q$ is exactly equal to the tangent space of this facet. Note that the
tangent space to a $k$ face is equal to the intersection of the
tangent spaces to all the maximal cells that contain it, so that
the statement extends to lower dimensional strata of the boundary of
$Q$.

Since the zero section is a Lagrangian submanifold of $T \bR^{n}$, the
polytope $Q$ may be thought of as a Lagrangian ball with boundary on
$M_{\infty}$; we denote this ball by $\cL_{\infty}$.
\end{rem}
Consider the Hamiltonian function
\[ H_{\infty}(u_1, \ldots , u_{i}) =  - \pi \sum_{i=1}^{n} {u_{i}}^2 \]
and let ${\phi^{1}_{\infty}}$ be its time--$1$ Hamiltonian flow. In the
universal cover, an explicit formula for ${\phi^{1}_{\infty}}$ is
given by
\[  {\phi^{1}_{\infty}}(u_1, \theta_1,   \ldots , u_n, \theta_n) = (u_1, \theta_1 - 2  \pi u_1, \ldots , u_n, \theta_n - 2  \pi u_n)\]
which we can write more conveniently as
\begin{equation}\label{twist} {\phi^{1}_{\infty}}(u , \theta) = (u ,\theta - 2  \pi u) . \end{equation}
\begin{lem}\label{infbound}
${\phi^{1}_{\infty}}(\partial \cL_{\infty} ) \subset M_{\infty} $
\end{lem}
\begin{proof}
Note that Equation \eqref{twist} implies that if $u$ is an integral point of $\bR^{n}$, then $\phi^{1}_{\infty}$ pointwise fixes the inverse image of $u$ under the projection map.  In particular, the inverse image of every vertex of $Q$ is fixed.

Consider a top dimensional cell of $Q$.  After a suitable translation by an element of $\bZ^{n}$, we may assume that one of its vertices is the origin.  But
the restriction of ${\phi^{1}_{\infty}}$ to the zero section of a
linear subspace is just the map
\[ (v,0) \to (v,-2 \pi v) .\]
So the subset of $M_{\infty}$ which lies over a top dimensional cell is preserved by  $\phi^{1}_{\infty}$.  Since  $\partial \cL_{\infty} \subset M_{\infty}$, this establishes the result for the restriction of $\partial \cL_{\infty}$ to the top dimensional cells.  The lemma follows from \fullref{Tsk}.
\end{proof}

We will denote ${\phi^{1}_{\infty}}(\cL_{\infty})$ by
$\cL_{\infty}(1)$.  The previous lemma in fact shows that the time $l$
flow satisfies
\[ {\phi^{l}_{\infty}} (\partial \cL_{\infty}) \subset M_{\infty} \]
for every integer $l$, so we have well defined Lagrangian balls
$\cL_{\infty}(l)$ with boundaries on $M_{\infty}$.

\subsection{A preliminary computation of Floer groups} \label{prelim}

Even though we don't have smooth boundary conditions yet, we will set
out to compute the Lagrangian Floer homology groups of the pairs
$(\cL_{\infty}(l_1), \cL_{\infty}(l_2))$.  We are missing the usual
compactness and transversality results that guarantee that our answers
will be invariant under Hamiltonian perturbations but we will proceed
regardless with the construction.  We will justify the use of small
Hamiltonian perturbations of the interior points by constructing admissible Lagrangians in Corollaries
\ref{prefloer1} and \ref{prefloer2}, which will also resolve the issue of boundary
intersection points.

Note that Lagrangian sections of the cotangent bundle have a natural
grading coming from the Morse index \acite{seidelGL}*{Example 2.10}. In
particular, whenever two such sections $\sL$, $\sL'$ intersect
transversely at a point $p$, we may write them locally as the graphs
of $df$ and $df'$ for smooth real valued functions $f$ and $f'$ such
that $p$ is a non-degenerate critical point of $f'-f$.  As an element
of $CF_{*}(\sL,\sL')$, the class of $p$ will have degree equal to its
Morse index as a critical point of $f'-f$.  We write
$\mu_{p}(\sL,\sL')$ for this index. These choices carry naturally to
the quotient of the cotangent bundle by a lattice coming from an
affine structure on the base \cite{leung} to give a canonical choice
of grading on sections of Lagrangian torus bundles.

We also need to resolve the issue of the status of intersection points
that occur on the boundary.  Our temporary ad-hoc prescription, which
we will justify in \fullref{comp}, is that boundary intersection
points are included in $CF_{*}(\cL_{\infty}(l_1), \cL_{\infty}(l_2))$
if and only if $l_1 < l_2$.  We will also stipulate that the group 
$CF_{*}(\cL_{\infty}(l), \cL_{\infty}(l))$ is isomorphic to $\bC$
concentrated in degree $n$.
\begin{lem}
The Floer homology groups between the Lagrangians $\cL_{\infty}(j)$
are given by
\[ HF_{n}(\cL_{\infty}(l_1),\cL_{\infty}(l_2)) = \bigoplus_{p \in Q \cap \frac{1}{l_2-l_1} \bZ^n} \bC \cdot [p]\]
with all other groups zero if $l_1 < l_2$, and
\[ HF_{0}(\cL_{\infty}(l_1),\cL_{\infty}(l_2)) = \bigoplus_{p \in ( Q - \partial Q) \cap \frac{1}{l_2-l_1} \bZ^n} \bC \cdot [p] \]
with all other groups zero if $l_1 > l_2$.
\end{lem}
\begin{proof}
First we reduce to the case where the first Lagrangian is the zero
section by applying an appropriate twist.  In other words, we have a
natural graded isomorphism
\[ CF_{*} (\cL_{\infty}(l_1), \cL_{\infty}(l_2)) \cong CF_{*} (\cL_{\infty}(l_1+i), \cL_{\infty}(l_2+i)) \]
for any integer $i$. Considering the case $i=-l_1$ reduces the
computation to one we've already done.

Indeed, we have already remarked in our proof of \fullref{infbound}
that the intersection points of $\cL_{\infty}$ and $\cL_{\infty}(l)$
correspond to the $\frac{1}{l} \bZ^{n}$ points of $Q$.  For
simplicity, we first discuss the case where $l>0$.  Note that the lift
of $\cL_{\infty}(l)$ to $T^{*}\bR^{n}$ which intersects the zero
section at a point $p = (p_1, \ldots, p_n)$ is given by the
differential of the function $- l\pi \sum_{i=1}^{n} (u_{i}-p_i)^{2}$.
Therefore, since the Morse index of this function at $p$ is $n$ if
$0<l$, the Floer complex is concentrated in degree $n$, where it is
given by
\[ CF_{n}(\cL_{\infty},\cL_{\infty}(l)) = \bigoplus_{p \in Q \cap \frac{1}{l} \bZ^n}  \bC \cdot [p].\]
The same computation yields that if $l<0$ the Floer complex is
concentrated in degree $0$, ie,
\[ CF_{0}(\cL_{\infty},\cL_{\infty}(l)) = \bigoplus_{p \in  (Q -\partial Q) \cap \frac{l}{l} \bZ^n} \bC \cdot [p] .\]
Since all these complexes are concentrated in one degree, the
differential is necessarily trivial.
\end{proof}

We will now pass to cohomology in order to compute the cup product.
Recall that the degree of a transverse intersection point $p \in \sL
\cap \sL'$ in cohomology is
\[ I_{p}(\sL,\sL')  = n - \mu_{p}(\sL,\sL').\]
We will compute the cup product
\[ HF^{*}(\cL_{\infty}(l_1),\cL_{\infty}(l_2)) \otimes HF^{*}(\cL_{\infty}(l_2),\cL_{\infty}(l_3))  \to HF^{*}(\cL_{\infty}(l_1),\cL_{\infty}(l_3)) ,\]
by counting holomorphic triangles connecting three intersection
points.  Assume the integers $l_1$, $l_2$, and $l_3$ are all distinct.

\begin{lem}
If $p=(p_1,\cdots, p_n)$ is an intersection point of
$\cL_{\infty}(l_1)$ and $\cL_{\infty}(l_2)$, and $q=(q_1, \cdots,
q_n)$ is an intersection point of $\cL_{\infty}(l_2)$ and
$\cL_{\infty}(l_3)$, then there exists at most one intersection point
$r$ of $\cL_{\infty}(l_1)$ and $\cL_{\infty}(l_3)$ such that
\[ \cM(p,q,r) \neq \emptyset .\]
Furthermore, if the intersection point $r$ exists, it is given by the
formula
\begin{equation} \label{add} r = \frac{(l_2-l_1) p + (l_3-l_2) q}{l_3-l_1}. \end{equation}
\end{lem}
\begin{rem}
The reader should note the similarity with Equation \eqref{fam}.
\end{rem}
\begin{proof}
There is a topological obstruction to the existence of a holomorphic
triangle with the appropriate boundary conditions which we now
describe.

After choosing a lift $\tilde{p}$ of $p$ to $T^{*} \bR^{n}$, there are
uniquely determined lifts of each Lagrangian.  In particular, the
lifts of $\cL_{\infty}(l_2)$ and $\cL_{\infty}(l_3)$ intersect at a
lift $\tilde{q}$ of $q$, while the lifts of $\cL_{\infty}(l_1)$ and
$\cL_{\infty}(l_3)$ intersect at most in a unique points which we call
$\tilde{r}$.  For the next few paragraphs, we will compute everything
in term of these lifts.

Using the metric to identify the cotangent fibre with the tangent
space, and the vector space structure on the base to identify each
tangent space with that of the origin, we may write $ \tilde{p} = (p,
\tilde{p}_{1}, \ldots, \tilde{p}_{n})$.  The lift of
$\cL_{\infty}(l_1)$ can be identified with the graph of the affine
transformation
\[ (x_1, \cdots , x_n) \to (-l_1(x_1-p_1) + \tilde{p}_{1}, \ldots,- l_1(x_n-p_n) + \tilde{p}_{n}) \]
and similarly for the lift of $\cL_{\infty}(l_2)$.  Since $q$ lies on
$\cL_{\infty}(l_2)$, its lift is therefore
\[ \tilde{q} = (q, -l_2(q_1-p_1) + \tilde{p}_{1}, \ldots, -l_2(q_n-p_n) + \tilde{p}_{n}). \]
This allows us to conclude that the lift of $\cL_{\infty}(l_3)$ that
we're considering is the graph of
\begin{align*}
x_i & \to -l_3(x_i-q_i) + \tilde{q}_{i} \\ & = -l_3(x_i-q_i) -l_2(q_i-p_i) +
\tilde{p}_{i}.\end{align*} Note that our sought after intersection point $r$ is
given by the solution to the system of linear equations
\[- l_1(r_i-p_i) + \tilde{p}_{i} = - l_3(r_i-q_i) - l_2(q_i-p_i) + \tilde{p}_{i} \, \, , \, \, 1 \leq i \leq n ,\]
which is clearly given by Equation \eqref{add}.
\end{proof}

\subsection{Cup product} \label{cup_product}

We must now count holomorphic triangles with appropriate boundary
conditions in order to compute the cup product.  As in the previous
section, any statements that are not justified here are handled in
future sections, in particular \fullref{comp}.  As we're working
with products of Lagrangians in different copies of $\bC^{\star}$, we
appeal to the fact that the holomorphic triangle admits no
deformations, so proving regularity for a pseudo-holomorphic map
\[ u\co (T, \partial T) \to (M, \cup_{i} L_i) \]
amounts to proving the surjectivity of the $\bar{\partial}$ operator
\[ W^{1,p}((T, \partial T),  (u^{*}TM, \cup_{i} u^{*}TL_i)) \to L^{p}(T,  \Omega^{0,1}(T) \otimes u^{*}TM ).\]
If $L_i \into M$ are given by products of Lagrangians in different
factors of the symplectic manifold $M$, then the above map splits into
direct summands, and surjectivity amounts to surjectivity for each of
the summands of $u$.

It is well known that non-constant holomorphic polygons on Riemann
surfaces are regular, though we will only give a proof in this
specific situation.
\begin{lem}
Consider three lines $L_1$, $L_2$ and $L_3$ in $\bR^{2} = \bC$ with
rational slope which intersect at three distinct points $(p,q,r)$.  If
there is a holomorphic triangle in $\cM(p,q,r)$, then it is
necessarily regular.
\end{lem}
\begin{proof}
Note that the usual index theoretic argument shows that such a
holomorphic triangle has index $0$.  In particular, it suffices to
show the operator
\[ W^{1,p}((T, \partial T),  (\bC, \cup_{i} TL_i)) \to L^{p}(T,  \bC )\]
is injective.

Assume that $F$ is an element of the kernel, and choose an integer $K$ such
that the images of the tangent lines $TL_i$ under the map $z \to
z^{K}$ is the real axis.  Note that $F^{K}$ will therefore be a
holomorphic map from $T$ to $\bC$ which takes $\partial T$ to $\bR$.
By the maximum principle, all such maps are constant.  However, since
$F$ has finite $W^{1,p}$ norm, and hence finite $L^p$ norm with respect to an infinite measure (because of the strip-like ends), it must therefore be identically $0$.
\end{proof}

\begin{lem}
Assume $p$ and $q$ are such that $r = \frac{(l_2-l_1) p + (l_3-l_2)
q}{l_3-l_1}$ lies in the polytope $Q$ and that the integers $l_1$, $l_2$, and $l_3$ are distinct.  If $l_1<l_2$ then
\[ \cM(p,q,r) \neq \emptyset \]
if and only if
\begin{equation}\label{cup1}  l_3  < l_1 < l_2  \, \, \mathrm{or} \, \, l_1 < l_2 < l_3  .\end{equation}
If $l_2<l_1$, then the moduli space is non-empty if and only if
\begin{equation}\label{cup2} l_2 < l_3 < l_1 .\end{equation} 
\end{lem}
\begin{proof}
By the preceding discussion, the existence of a holomorphic triangle
in the total space with boundary on our given affine subspaces is
equivalent to the existence of $n$ holomorphic triangles with boundary
on a configuration of straight lines in $\bR^{2}$.  But Riemann's
mapping theorem says that such a triangle exists if and only if
orientations are preserved.  Moreover, since in each factor such a
triangle is unique if it exists, there can be at most one holomorphic
disc with appropriate boundary conditions.  One might also be
concerned that the holomorphic triangles we're finding may go outside
the the inverse image of $Q$ even when all the intersection points between
the different Lagrangians lie within $Q$.  But this doesn't happen
because $Q$ is convex. Indeed, representing the triangle with three
marked points as a geodesic triangle $T$ in the plane, we see that
each component $u_i$ of the holomorphic map $u\co T \to \bC^{n}$ is in
fact given by an affine transformation of the plane.  In particular,
the image of $T$ is the flat triangle determined by the points $(
\tilde{p}, \tilde{q}, \tilde{r} )$, hence projects to the flat
triangle determined by $(p,q,r)$.  Convexity implies that this
triangle is contained in $Q$.

\begin{figure}\small
\begin{center}
$\begin{array}{c@{\hspace{.5in}}c} 
\begin{picture}(0,0)%
\includegraphics{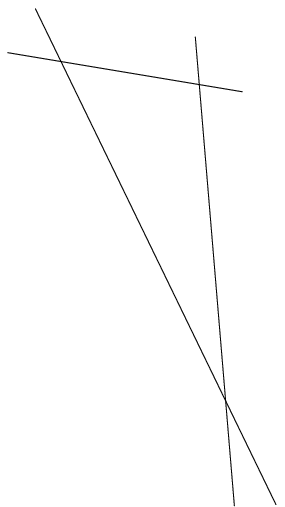}%
\end{picture}%
\setlength{\unitlength}{1973sp}%
\begin{picture}(2661,4795)(4051,-13348)
\put(4051,-8836){\makebox(0,0)[lb]{\smash{{\SetFigFont{6}{7.2}{\familydefault}{\mddefault}{\updefault}{\color[rgb]{0,0,0}$\tilde{p}$}%
}}}}
\put(6076,-10411){\makebox(0,0)[lb]{\smash{{\SetFigFont{6}{7.2}{\familydefault}{\mddefault}{\updefault}{\color[rgb]{0,0,0}$\cL_{\infty}(l_3)$}%
}}}}
\put(6001,-9211){\makebox(0,0)[lb]{\smash{{\SetFigFont{6}{7.2}{\familydefault}{\mddefault}{\updefault}{\color[rgb]{0,0,0}$\tilde{r}$}%
}}}}
\put(4801,-9061){\makebox(0,0)[lb]{\smash{{\SetFigFont{6}{7.2}{\familydefault}{\mddefault}{\updefault}{\color[rgb]{0,0,0}$\cL_{\infty}(l_1)$}%
}}}}
\put(4200,-10411){\makebox(0,0)[lb]{\smash{{\SetFigFont{6}{7.2}{\familydefault}{\mddefault}{\updefault}{\color[rgb]{0,0,0}$\cL_{\infty}(l_2)$}%
}}}}
\put(6301,-13086){\makebox(0,0)[lb]{\smash{{\SetFigFont{6}{7.2}{\familydefault}{\mddefault}{\updefault}{\color[rgb]{0,0,0}$\tilde{q}$}%
}}}}
\end{picture}%
&
\begin{picture}(0,0)%
\includegraphics{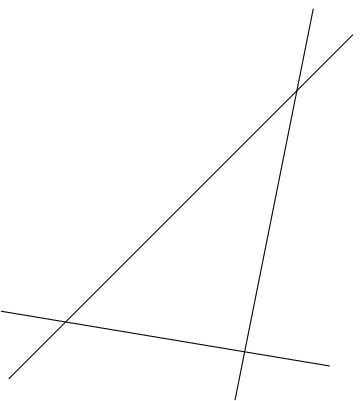}%
\end{picture}%
\setlength{\unitlength}{1973sp}%
\begingroup\makeatletter\ifx\SetFigFont\undefined%
\gdef\SetFigFont#1#2#3#4#5{%
  \reset@font\fontsize{#1}{#2pt}%
  \fontfamily{#3}\fontseries{#4}\fontshape{#5}%
  \selectfont}%
\fi\endgroup%
\begin{picture}(3399,3774)(4114,-9847)
\put(4126,-9286){\makebox(0,0)[lb]{\smash{{\SetFigFont{6}{7.2}{\familydefault}{\mddefault}{\updefault}{\color[rgb]{0,0,0}$\tilde{p}$}%
}}}}
\put(6676,-6736){\makebox(0,0)[lb]{\smash{{\SetFigFont{6}{7.2}{\familydefault}{\mddefault}{\updefault}{\color[rgb]{0,0,0}$\tilde{r}$}%
}}}}
\put(6526,-9761){\makebox(0,0)[lb]{\smash{{\SetFigFont{6}{7.2}{\familydefault}{\mddefault}{\updefault}{\color[rgb]{0,0,0}$\tilde{q}$}%
}}}}
\put(6676,-8536){\makebox(0,0)[lb]{\smash{{\SetFigFont{6}{7.2}{\familydefault}{\mddefault}{\updefault}{\color[rgb]{0,0,0}$\cL_{\infty}(l_2)$}%
}}}}
\put(5026,-7736){\makebox(0,0)[lb]{\smash{{\SetFigFont{6}{7.2}{\familydefault}{\mddefault}{\updefault}{\color[rgb]{0,0,0}$\cL_{\infty}(l_3)$}%
}}}}
\put(4876,-9636){\makebox(0,0)[lb]{\smash{{\SetFigFont{6}{7.2}{\familydefault}{\mddefault}{\updefault}{\color[rgb]{0,0,0}$\cL_{\infty}(l_1)$}%
}}}}
\end{picture}%
\\ \mbox{ (1)} & \mbox{ (2)}
\end{array}$
\end{center}
\caption{Holomorphic discs with $l_1 < l_2 < l_3$ (1), and $ l_2 < l_3
< l_1$ (2)}
\label{holo_triangles}
\end{figure}

We have therefore reduced the problem to establishing whether the
triples $(\xi, \xi^{2}, 1)$ and $( (p_{i}, \tilde{p}_{i}), (q_{i},
\tilde{q}_{i}), (r_{i}, \tilde{r}_{i}) )$ have the same orientation,
assuming the latter consists of distinct points. It follows from
Formula \eqref{add} (see also, \fullref{holo_triangles}) that if
$l_1<l_2$, the orientations are the same if and only if
\[  l_3  < l_1 < l_2  \, \, \mathrm{or} \, \, l_1 < l_2 < l_3 \]
while if $l_2<l_1$, then the orientations are the same if and only if
\[ l_2 < l_3 < l_1 . \]
In particular, this condition is independent of which factor we're
considering, so the above is the necessary and sufficient condition
for the existence of a holomorphic disc with marked points on
$(\tilde{p}, \tilde{q}, \tilde{r})$ if these points are distinct.

We must now deal with the case where the marked points are not
distinct.  Observe that since all the slopes $l_i$ are assumed to be different, this can only occur if all three points are
equal.  We must appeal to a yet unproven invariance under Hamiltonian
deformations to interpret the product among Lagrangian submanifolds in
non-generic positions in terms of the product for nearby
Lagrangians. Note the Conditions \eqref{cup1} or \eqref{cup2} are clearly
invariant under small Hamiltonian deformations.  In particular, given
a triple intersection point between $\cL_{\infty}(l_1)$,
$\cL_{\infty}(l_2)$, and $\cL_{\infty}(l_3)$, we may locally deform
one of the Lagrangians, through Lagrangians which are also given as
$\cL_{\infty}(l)$ for some real number $l$, preserving Conditions \eqref{cup1} or
\eqref{cup2}.  But the argument we gave when $l_i$ are distinct integers
did not depend on their integrality, so we reach the desired
conclusion that Conditions \eqref{cup1} or \eqref{cup2} suffice to
describe the existence of a regular holomorphic triangle.
\end{proof}
\begin{rem}
The perturbation argument that we're using breaks down whenever $p=q=r$ is a point on the boundary.  In \fullref{admiss-sect}, we will ``push intersection points to the interior,'' and thereby legitimate the perturbation argument.
\end{rem}

If two of the Lagrangian sections are equal, recall that we have made
an as yet unjustified claim that the corresponding Floer cohomology
group is concentrated in degree $0$ and is isomorphic to a copy of
$\bC$.  We will justify in \fullref{comp} that under the cup
product, this group acts by scalar multiplication on the Floer
cohomology group of a pair of Lagrangian manifolds. 

To complete our analysis of the presence of holomorphic triangles, we
must find appropriate conditions for $ \tilde{r}$ to lie within the
prescribed region.  We will only do this for a special case.
\begin{lem} \label{prelim_product}
Assume that $l_1 \leq l_2 \leq l_3$.  The product in Floer cohomology
is given by
\begin{align}
 HF^{0}(\cL_{\infty}(l_1), \cL_{\infty}(l_2)) \otimes
HF^{0}(\cL_{\infty}(l_2), \cL_{\infty}(l_3)) &\to
HF^{0}(\cL_{\infty}(l_1), \cL_{\infty}(l_3)) \nonumber\\
 \label{product} [ p ] \otimes [ q ] &\mapsto \pm \left[ \frac{(l_2-l_1) p + (l_3-l_2) q}{l_3-l_1} \right].
\end{align}
\end{lem}  
\begin{proof}
Since $Q$ is convex, and all the coefficients in Equation \eqref{add}
are positive, we conclude that $\tilde{r}$ projects to an intersection
point of $\cL_{\infty}(l_1)$ and $\cL_{\infty}(l_3)$.  In the case of
strict inequalities among the $l_j$, Equation \eqref{cup1} implies that
there is necessarily a holomorphic triangle connecting $p$, $q$ and
$r$.  The case where two of the Lagrangians may be equal follows from
the previous discussion.
\end{proof}

We also know from Equation \eqref{cup2} that there might be other
non-trivial cup products if the condition $l_1 \leq l_2 \leq l_3$ does
not hold. The point is that these are the same products that show up
on the algebraic side. We will take this up again in \fullref{end}.

We will now prove that we may choose all signs in Equation
\eqref{product} to be positive. 
\begin{lem}
There is a choice of generators for $CF^{*}(\cL(l_i), \cL(l_j))$ such
that all the signs that arise in the cup product are positive.
\end{lem}
\begin{proof}
We already know that
\[ [p] \otimes [q] \mapsto  \pm [r] \]
and it remains to prove that all signs may be chosen to be positive.
However, since any configuration of three Lagrangian sections with
$l_1 < l_2 <l_3$ is isotopic through Lagrangian sections to any other
configuration and the cup product is invariant under such Lagrangian
isotopies, we conclude that all signs are either positive or negative.
If necessary, we complete the argument by replacing every generator by its negative.
\end{proof}

In fact, this argument extends to products in which the triple of
Lagrangians does not necessarily satisfy the condition $l_1 < l_2 <
l_3$.  Indeed, one can then simply dualize the operation that induces
the cup product to reduce every other possible configuration to the
case $l_1 <l_2 <l_3$.  For example, if $l_3 < l_1 < l_2$, then, restricting to the interior intersection points, the map
\[ CF^{*}(\cL(l_1), \cL(l_2)) \otimes CF^{*}(\cL(l_2), \cL(l_3)) \to CF^{*}(\cL(l_1), \cL(l_3)) \] 
can be obtained by dualizing the map
\[ CF^{*}(\cL(l_3), \cL(l_1)) \otimes CF^{*}(\cL(l_1), \cL(l_2)) \to CF^{*}(\cL(l_3), \cL(l_2)).\] 
But we know from the previous lemma that with our choice of generators, the sign that appears in this last product is positive,
hence so it is for the first one.  Once this is established for the interior intersection points, it immediately follows that it is true when boundary intersection points are involved as well because the signs are determined by local considerations which are insensitive to the distinction between the interior and the boundary.

\section{Tropical localization}

Mikhalkin used \fullref{GH} to study the topology of smooth
complex hypersurfaces $M$ of $\ctorus{n}$ by producing a subset of
$\bR^{n}$ which interpolates between the polyhedral complex $\Pi$
(thought of as the tropical amoeba of a non-Archimedean variety) and
the $\Log$--amoeba of $M$.  He observed that the tropical amoeba has a
complex analogue to which the hypersurface $M$ degenerates.  We will
give a slight modification of his proof which will allow us to obtain
a stronger result about the symplectic structures.

\subsection{Construction}

Recall from \fullref{TG} that each component of the complement of
$\log(t) \Pi = \Pi_{t}$ corresponds to a vertex $\alpha$ of a coherent
subdivision.  We defined $C_{\alpha}$ in terms of the linear function
attached to $\alpha$, but one might equivalently define $\log(t) C_{\alpha} = C_{\alpha,t}$ as the image of the open subset of $\ctorus{n}$ where the monomial
$t^{-\nu(\alpha)} z^{\alpha}$ dominates all others.  For the purposes of this section, we will need some bounds on polyhedral geometry of
$P$.  For simplicity, we will assume that $c_{\alpha}=1$.

Let $N$ be the maximum of the $\ell^1$ norms of all vectors $\alpha - \beta$
for $\alpha$, $\beta$ neighbouring vertices of the subdivision, and of all vectors $\alpha \in A$.  Note that $N$ is also an upper bound for the Euclidean norm of such vectors.  Assume that $\nu$
induces a triangulation, and that the affine map from the standard
$k$--dimensional vector space to $\bR^{n}$ induced by every simplex of
the subdivision of $A$ distorts length by an amount bounded by $\rho >1$.  Note that such a bound exists simply because these affine maps are all injective.

If $\alpha$ and $\beta$ are distinct elements of $A$, let $\cH(\alpha, \beta)$ denote the hyperplane defined by the equation
\[ \langle \alpha, u \rangle - \nu(\alpha) = \langle \beta, u \rangle - \nu(\beta) .\]
Note that every face of $\Pi$ is supported by such a hyperplane. 
\begin{lem} \label{distance_hyperplane}
There exists a constant $c>0$ such that for $\epsilon$ sufficiently small
\[ d(p,C_{\alpha}) \geq \epsilon \Rightarrow d(p, \cH(\alpha, \beta)) \geq 2 c \epsilon \]
for every $p \in C_{\beta}$ and for all pairs $\alpha \neq \beta$.
\end{lem}
\begin{proof}
Since $A$ only has finitely many elements, it suffices to found a bound for each of them and then take the minimum. Accordingly, fix a component $C_{\beta}$.

Choose a constant $\delta$ such that $2 \delta$ neighbourhood of $C_{\alpha}$ can only intersect $C_{\beta}$ if their boundaries share a face.  Since $\Pi$ is a finite cell complex and $\cH(\alpha, \beta)$ cannot asymptotically approach $C_{\beta}$ at infinity, there exists a constant $K$ such that $d(p,\cH(\alpha, \beta)) \geq K$ whenever $p$ does not lie in the $\delta$ neighbourhood of $C_{\alpha}$.

Let us assume for simplicity that the face $\sigma$ shared by $C_\alpha$ and $C_\beta$ is bounded.  In this case, the boundary of the $\delta$ neighbourhood of $\sigma$ intersects $C_{\beta}$ in a compact subset which we will denote $S$.  Note that convexity of $C_{\beta}$ implies that $C_{\beta} \cap \cH(\alpha, \beta) = \sigma$.  In particular, the distances to $\cH(\alpha,\beta)$ and $\sigma$ are bounded above and below on $S$ by non-zero constants, so an appropriate ratio yields the desired constant $c$.  Consider $p \in C_\beta - C_\alpha$ such that $d(p,\sigma) = \epsilon < \delta$.  There are points $q$ on $\sigma$ and $r$ on $\cH(\alpha, \beta)$ which realize the distance from $p$ to these respective sets.  The ratio between these distances can be computed from the angles of the triangle with corners $p$, $q$ and $r$.  By extending the segment from $p$ to $q$, we eventually reach a point $p'$ on $S$.  Note that $q$ is the point of $\sigma$ closest to $p'$.  If $r'$ is the orthogonal projection of $p'$ on $\cH(\alpha, \beta)$, then the right triangles with corners $(p,q,r)$ and $(p',q,r')$ are similar.   In particular, the ratios of their sides are equal.  Since the ratio of the sides of the triangle $(p',q,r')$ is bounded by $2c$, the same bound works for the triangle $(p,q,r)$.  Further, since 
\[d(p,q) \geq d(p, C_{\alpha}) = \epsilon,\]
we conclude that 
\[ d(p, \cH(\alpha, \beta)) = d(p,r) \geq 2 c d(p,q) \geq 2 c \epsilon. \]
If the intersection between $C_{\alpha}$ and $C_\beta$ is not compact, we choose a ball $B$ that for every point $q_1$ of $(C_{\alpha} \cap C_{\beta})- B$ there is a point $q_2$ in $C_{\alpha} \cap C_{\beta} \cap B$ such that there is an isometry of the $2 \delta$ neighbourhoods of $q_1$ and $q_2$ which takes the intersection of $\Pi$ with one neighbourhood to its intersection with the other.  The existence of such a ball is a consequence of the fact that $\Pi$ is a finite polyhedral complex.  If $p$ lies in the $\delta$ neighbourhood of $q$, then the point of $\cH(\alpha, \beta)$ nearest to $p$ lies in the $2 \delta$ neighbourhood of $q$.  So the problem is entirely local, and the bound obtained by considering the intersection of $B$ with the boundary of the $\delta$ neighbourhood of $C_{\alpha}$ in $C_{\beta}$ yields a bound which is also valid for the complement of $B$.  The rest of the argument proceeds as in the compact case.
\end{proof}

Now pick $\epsilon$ such that $\delta > \epsilon >0$.  By \fullref{GH}, there exists a $T$ such
that $d(\cA_t, \Pi_t) < \log(t) \epsilon$ for all $t>T$.  Choose such
$\epsilon$ and $T>1$ such that the following additional conditions are
satisfied for all $t>T$
\begin{align}
\frac{e^{-c \epsilon
\log(t)}}{\epsilon \log(t)} & <  \frac{1}{40\left|A\right| \rho}
\label{e-2}\\ e^{-c\epsilon \log(t)} & < 
\frac{1}{5\left|A\right|^{2} \rho N}. \label{e-3}
\end{align}
In addition, choose nowhere negative $C^{\infty}$ functions $\phi_{\alpha}$ on
$\bR^{n}$ such that the following properties hold
\begin{align}
d(p,C_{\alpha,t}) \leq \frac{\epsilon \log(t)}{2} & \Leftrightarrow  \phi_{\alpha}(p) = 0
\label{phi-1} \\ d(p, C_{\alpha,t}) \geq \epsilon \log(t) & \Leftrightarrow 
\phi_{\alpha}(p) = 1 \label{phi-2} 
\end{align}
\begin{equation} \sum_{i=1}^{n} \left| \frac{\dd
\phi_\alpha (p)}{\dd u_{i}} \right| <\frac{4}{\epsilon \log(t)}.
\label{phi-3} \end{equation}
We will also
abuse notation and write $\phi_{\alpha}(z)$ for $\phi_{\alpha}(\Log(z))$.  Further, we will assume that $\phi_{\alpha}$ and the norm of its derivatives are $C^{0}$ close to functions which depend only on the distance to $C_{\alpha}$.

We now consider the family of maps \[f_{t,s}=\sum_{\alpha \in A}
t^{-\nu(\alpha)} (1 - s \phi_{\alpha} (z) )z^{\alpha}.\] Our goal is
to prove the following:
\begin{prop} \label{family}
For large enough $t$, $f_{t,s}^{-1}(0)=M_{t,s}$ is a family of
symplectic hypersurfaces parametrized by $s$.
\end{prop} 

\begin{figure}
\begin{center}
$\begin{array}{c@{\hspace{.1in}}c@{\hspace{.1in}}c} 
\includegraphics[width=1.5in]{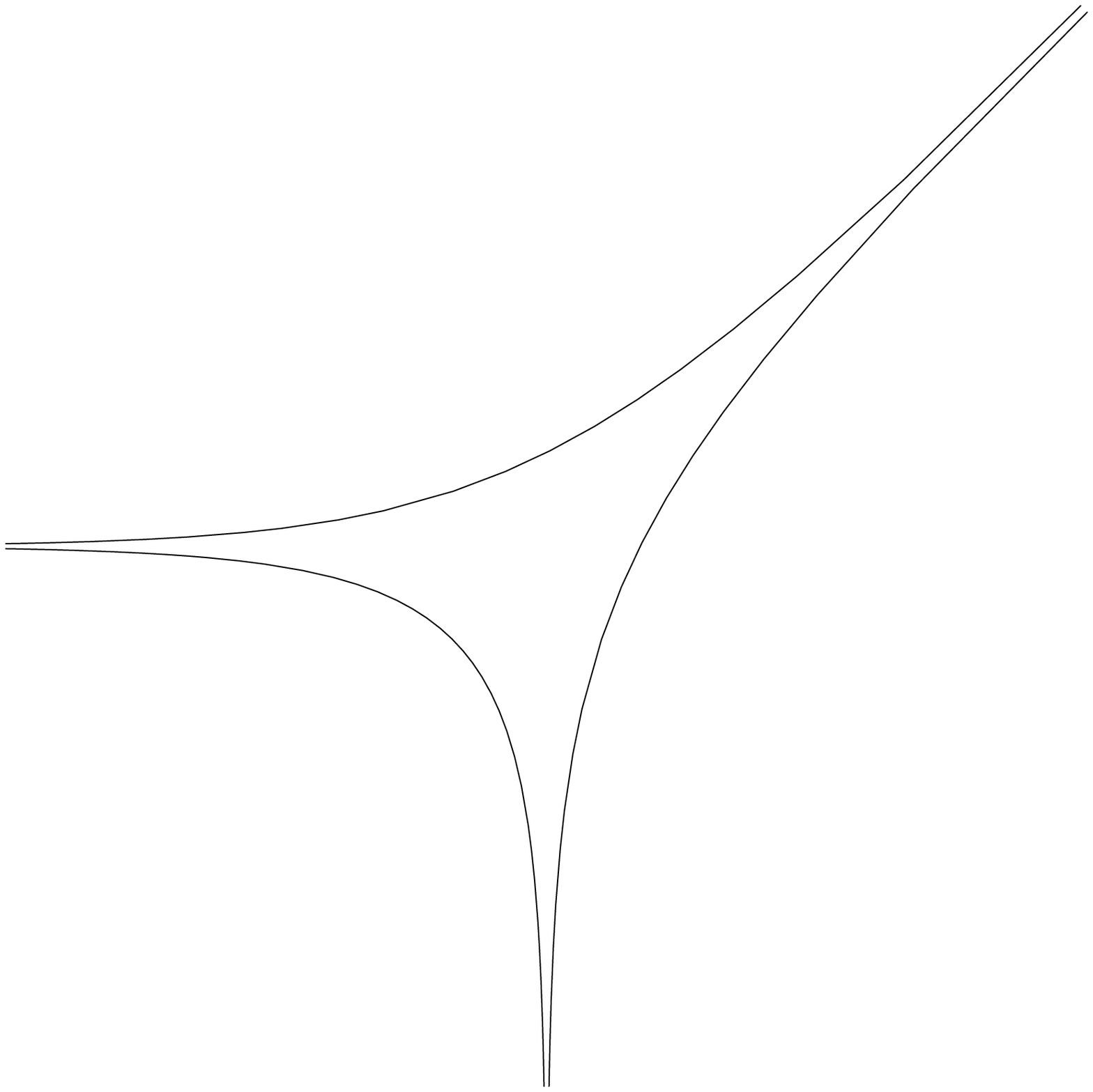} 
& 
\includegraphics[width=1.5in]{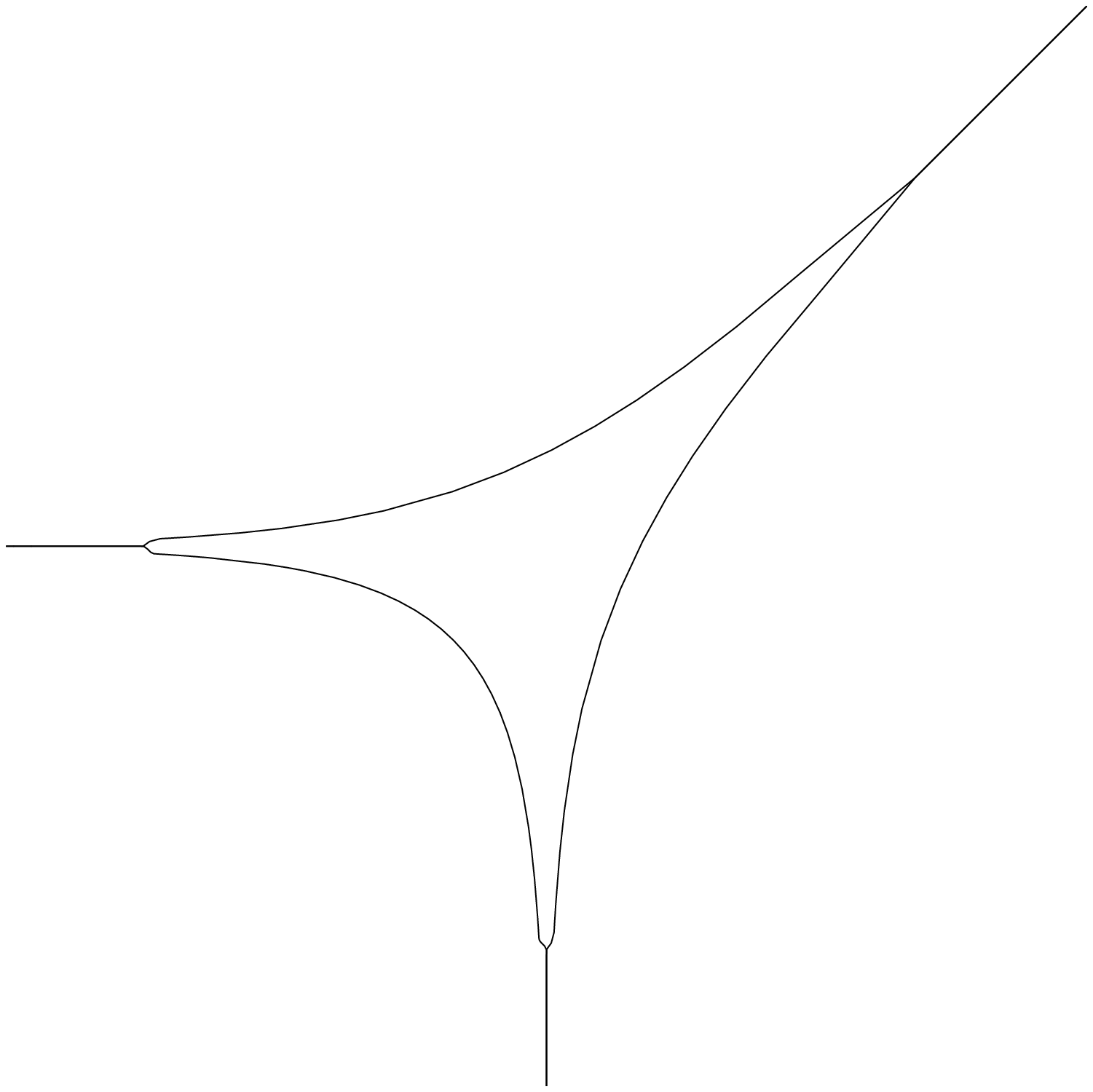} 
& 
\includegraphics[width=1.5in]{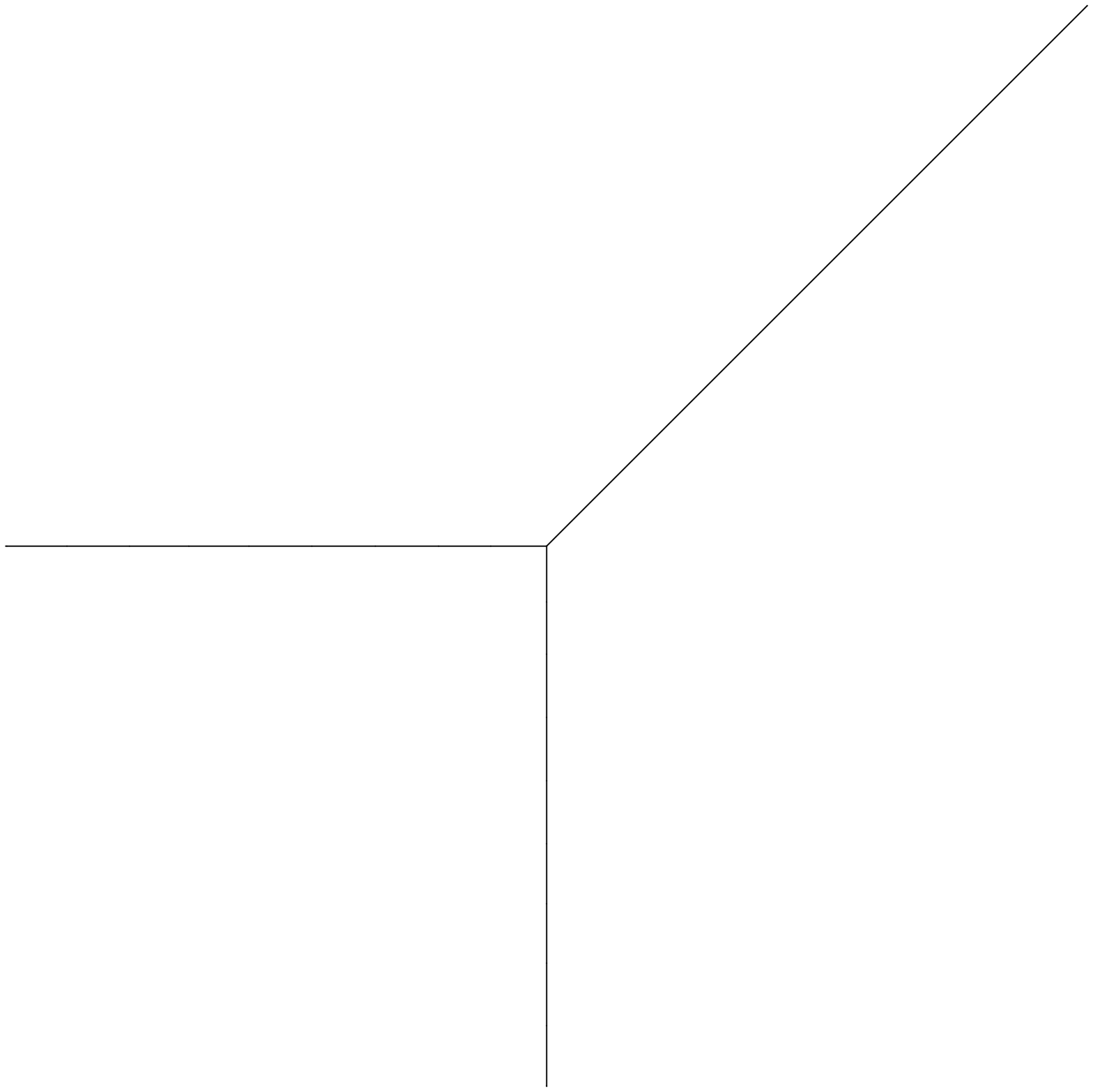}
\end{array}$
\end{center}
\caption{The amoeba of the standard hyperplane in $\ctorus{2}$, its
``tropical localization," and its tropical amoeba}
\label{amoebas-fig}
\end{figure}

We will call $f_{t,1}$ a ``tropical localization'' of the Laurent
polynomial $f$.  The choice of terminology should be clear from \fullref{amoebas-fig}.

It follows from the definition of $M_{t,1}$ that near a face $\sigma_t$
of $\Pi_t$, the equations that define $M$ depend only on those
coordinates which are orthogonal to $\sigma_t$.  In order to be more
precise, we consider a fixed cover of $\bR^{n}$.
\begin{defin}
Given any polytope $\tau$ of the subdivision of $P$, define $O_{\check{\tau}}$ to be the set of points $p \in \bR^n$ such that $\phi_{\alpha}(p) \neq 1$ for all $\alpha \in \tau$, and $\phi_{\alpha}(p) = 1$ for all $\alpha \notin \tau$.  
\end{defin}
Our choices of functions $\phi_{\alpha}$ ensure that 
\[ \bigcup_{\tau} O_{\check{\tau}} = \bR^n. \]
One way to see this is to observe that, although it is not an open set, we can still think of each $O_{\check{\tau}}$ as a neighbourhood of an open subset of  $\check{\tau_{t}} = \log(t) \check{\tau}$, the face of $\Pi_t$ dual to $\tau$.  Indeed, it is not hard to show that a point $p$ that is distance more than $\epsilon \log(t)$ away from $\check{\tau_{t}}$ must be distance at least $\epsilon \log(t)$ from some $C_{\alpha,t}$ with $\alpha \in \tau$, so $\phi_{\alpha}(p)=1$ by construction.  This proves that $O_{\check{\tau}}$ is contained in an $\epsilon \log(t)$ neighbourhood of $\check{\tau_{t}}$.  

However, since every face of $\partial \check{\tau}$ is dual to a polytope $\tau' \subset Q$ such that $\tau \subset \partial \tau'$, we know that near every face $\check{\tau}'_{t}$ of the boundary of $\check{\tau_t}$, there exists at least one vertex $\beta \notin \tau$ such that $\phi_{\beta}$ is not equal to $1$.  This implies that $\partial {\check{\tau_t}}$ is not contained in $O_{\check{\tau}}$.  But if we restrict to the complement of an $\epsilon \log(t)$ neighbourhood of $\partial \check{\tau_t}$ in $\check{\tau_t}$, then $\phi_{\beta}=1$ for every $\beta \notin \tau$.

In addition, note that since we're working with maximal subdivisions, every polytope $\tau$ of the subdivision of $P$ induces a product decomposition
\[ \ctorus{n} =  \ctorus{\tau} \times \ctorus{\check{\tau}} \]
where $\check{\tau}$ is the face of $\Pi$ which is dual to $\tau$.
We therefore have the following description.
\begin{lem} \label{semi-local}
In $O_{\check{\tau}}$, $M_{t,1}$ is a product
\[ M_{t,1} = M^{\tau}_{t,1} \times U ,\] with $ M^{\tau}_{t,1}$ a hypersurface in $\ctorus{\tau}$ and $U$ a subset of  $ \ctorus{{\check{\tau}}}$.  The set $U$ may itself be decomposed at $\bT^{\check{\tau}} \times U'$ with $\bT^{{\check{\tau}}}$ the real subtorus of $\ctorus{\check{\tau}}$.  \noproof
\end{lem}

Away from its boundary, $U'$ is an open set of $(\bR^{+})^{k}$ and may
be identified with the open subset of $\check{\tau}$ consisting of points
which lie sufficiently far away from the boundary. Furthermore, in the identification of $\ctorus{n}$ with $T
\bR^{n} / \bZ^{n}$, the torus $\bT^{\check{\tau}}$ corresponds to the
tangent space of $\check{\tau}$.

To prove \fullref{family}, we will use a simple
observation which was used by Donaldson in \cite{donaldson}: In order
to show that $f^{-1}(0)$ is symplectic it suffices to prove that on
the preimage of $0$, $\left|\overline{\dd}f\right| < \left|
\dd f \right|$, where the norms are computed with respect to the
unique metric determined by the complex and symplectic structures.  We
compute that
\[ \dd f_{t,s} = \sum_{\alpha \in A}  t^{-\nu(\alpha)} \dd z^{\alpha}  - s \sum_{\alpha \in A} t^{-\nu(\alpha)}\left( \phi_{\alpha}(z) \dd z^{\alpha} + z^\alpha \dd\phi_\alpha(z) \right) \]
while
\[ \overline{\dd}f_{t,s} = s \sum_{\alpha \in A}  t^{-\nu(\alpha)} z^\alpha \overline{\dd}\phi_\alpha(z). \]
Since our norm comes from a hermitian metric, and $\phi_\alpha$ is
only a function of the norm of $z$, it suffices to show that
\begin{equation} \label{symplectic_ineq} \left|\sum_{\alpha \in A}  t^{-\nu(\alpha)} \dd z^{\alpha} \right| > s \left(  \left|\sum_{\alpha \in A}  t^{-\nu(\alpha)} \phi_{\alpha}(z) \dd z^{\alpha}\right| + 2 \left| \sum_{\alpha \in A}  t^{-\nu(\alpha)}  z^\alpha \dd\phi_\alpha(z) \right| \right) \end{equation}
hence it suffices to show the above inequality for $s=1$.

The idea behind our choice of cut-off functions $\phi_{\alpha}$ is that $\phi_{\alpha}$ is only non-zero whenever the corresponding monomial contributes negligibly to $f_{t,s}$.  Formally we have the following:
\begin{lem} \label{expo_decay}
If $\Log(z)=p$ is in $C_{\beta,t}$ and  $\phi_{\alpha}(p) \neq 0$ then 
\[ \left| \frac{t^{-\nu(\alpha)} z^{\alpha}}{t^{-\nu(\beta)} z^{\beta}} \right| < e^{- c \epsilon \log(t) \left|\alpha-\beta\right|} .\]
\end{lem}
\begin{proof}
Consider the shortest segment with one endpoint at $p$ and the other at $q$, the closest point in $\cH(\alpha, \beta)$.  We have
\[ \langle q , \beta \rangle - \log(t) \nu(\beta) =  \langle q , \alpha \rangle - \log(t) \nu(\alpha). \]
We compute
\begin{align*}
\left| \frac{t^{-\nu(\alpha)} z^{\alpha}}{t^{-\nu(\beta)} z^{\beta}} \right| & = 
\frac{t^{-\nu(\alpha)} \left| e^{\langle p, \alpha \rangle} \right| }{t^{-\nu(\beta)} \left|  e^{\langle p, \beta \rangle } \right| } \\
& = e^{\langle p, \alpha - \beta \rangle + \log(t) (-\nu(\alpha) + \nu(\beta))}\\
& =  e^{\langle q, \alpha - \beta \rangle + \log(t) (-\nu(\alpha) + \nu(\beta))} e^{\langle p-q, \alpha - \beta \rangle} \\
& =  e^{\langle p-q, \alpha - \beta \rangle} .
\end{align*}
Since the segment from $p$ to $q$ is normal to $\sigma_t$, the vector from $p-q$ is parallel to $\alpha - \beta$, so 
\begin{align*}
\langle p-q, \alpha - \beta \rangle & = - |p-q| |\alpha - \beta| \\
& \leq - c \epsilon \log(t)  |\alpha - \beta|
\end{align*}
where the last inequality follows from $\phi_{\alpha}(p) \neq 0$, and \fullref{distance_hyperplane}.
\end{proof}

We will only be considering points at which $f_{t,s}$ vanishes. The
next two lemmata establish constraints that such points must
satisfy.
\begin{lem}\label{close-amoeba}
If $\cA_{t,s}$ is the amoeba of the hypersurface $M_{t,s}$, then
$d(\cA_{t,s}, \Pi_t) < \log(t) \epsilon$ for all $s$.
\end{lem}
\begin{proof}
Note that if $d(p,\Pi_t) \geq \log(t) \epsilon$, then there exists a unique
monomial $t^{-\nu(\beta)} z^{\beta}$ such that $\phi_{\beta}=0$.  By \fullref{expo_decay} and Condition \eqref{e-3}, we have
\begin{align*}
\left| (1-s\phi_{\alpha}(z)) t^{-\nu(\alpha)} z^{\alpha}\right| & \leq \left|  t^{-\nu(\alpha)} z^{\alpha} \right|
\\ & < \frac{1}{5\left|A\right|^{2} \rho N}\left|t^{-\nu(\beta)} z^{\beta}\right| \\
& < \frac{1}{\left|A\right|}\left|t^{-\nu(\beta)}  z^{\beta}\right|
\end{align*}
for every $\alpha \neq \beta$.  Since there are exactly $\left|A\right|$ such terms in our expression for $f_{t,s}$ for which $\alpha \neq \beta$, we conclude that $f_{t,s}(z)$ cannot vanish.
\end{proof}
\begin{lem} \label{bound_f_t}
If $s$ and $z$ are such that $f_{t,s}(z)=0$, and $\Log(z) \in C_{\beta}$, then 
\[\left|f_{t,r}(z)\right| < \frac{\left| t^{-\nu(\beta)}z^{\beta} \right|}{5|A|\rho N}\] for
every $0 \leq r \leq 1$.  The same bound holds for:
\[ \left| \sum_{\phi_{\alpha}(z) \neq 1} t^{-\nu(\alpha)} z^{\alpha} \right| \]
\end{lem}
\begin{proof}
It suffices to observe that:
\begin{align*}
\left|f_{t,r}(z)\right| & = \left|f_{t,r}(z) - f_{t,s}(z)\right| \\
& = \left|(r-s) \sum_{\phi_{\alpha}(z) \neq 0} t^{-\nu(\alpha)}
\phi_{\alpha}(z) z^{\alpha}\right| \\ & \leq \sum_{\phi_{\alpha}(z) \neq 0 } \left| t^{-\nu(\alpha)} z^{\alpha}\right| \\ & \leq
\left|A\right| e^{-c\epsilon \log(t)} \left|t^{-\nu(\beta)}
z^{\beta}\right| \\ & \leq \frac{\left| t^{-\nu(\beta)}z^{\beta}
\right|}{5\rho \left|A \right| N} \end{align*}
The same computation yields the second part of the lemma:
$$
\eqalignbot{
\left| \sum_{\phi_{\alpha}(z) \neq 1} t^{-\nu(\alpha)} z^{\alpha}\right|  & = \left| s \sum_{0 <\phi_{\alpha}(z) <1}  \phi_{\alpha}(z) t^{-\nu(\alpha)} z^{\alpha} - \sum_{\phi_{\alpha}(z) = 1} t^{-\nu(\alpha)} (1-s\phi_{\alpha}(z)) z^{\alpha}  \right| 
\cr & \leq \sum_{\phi_{\alpha}(z) \neq 0 } \left| t^{-\nu(\alpha)} z^{\alpha}\right| 
\cr & \leq \frac{\left| t^{-\nu(\beta)}z^{\beta} \right|}{5\rho \left|A \right| N}} \proved$$
\end{proof}
We now complete the proof of this section's main result.
\begin{proof}[Proof of  \fullref{family}]
Assume that $f_{t,s}(z)=0$, and $\Log(z) \in O_{\check{\tau}}$.  By \fullref{close-amoeba}, $\tau$ consists of at least two vectors. Since the subdivision induced by $\nu$ is a triangulation, each subpolytope is in fact a simplex. We choose $\delta$ in $\tau$ such that $\Log(z) \in C_{\delta,t}$, which implies that $\left|t^{-\nu(\delta)} z^{\delta} \right| \geq \left| t^{-\nu(\beta)}z^{\beta} \right| $ for all $\beta$.  In addition, choose a vertex $\gamma $ in $\tau$ which differs from $\delta$, and re-write
\[ \sum_{\alpha \in A}  t^{-\nu(\alpha)} \dd z^{\alpha} = \sum_{\beta \in \tau}  t^{-\nu(\beta)}  z^{\gamma} \dd z^{\beta-\gamma} + \sum_{\beta \in \tau}  t^{-\nu(\beta)} z^{\beta-\gamma} \dd{z^{\gamma}} + \sum_{\alpha \in A -\tau}  t^{-\nu(\alpha)} \dd z^{\alpha} .\]
In order to establish Equation \eqref{symplectic_ineq}, we decompose its left hand-side as above, so it will suffice to prove the inequality:
\begin{align*}
\left|z^{\gamma} \sum_{\beta \in \tau} t^{-\nu(\beta)} \dd
z^{\beta-\gamma}\right| & >  \left|\sum_{\alpha \in A} 
t^{-\nu(\alpha)} \phi_{\alpha}(z) \dd z^{\alpha} \right| 
 +  2
\left| \sum_{\alpha \in A} t^{-\nu(\alpha)} \dd \phi_\alpha(z)
z^\alpha \right| \\ &  + \left| \sum_{\beta \in \tau} t^{-\nu(\beta)}
z^{\beta-\gamma} \dd{z^{\gamma}}\right| + \left| \sum_{\alpha \notin \tau}  t^{-\nu(\alpha)} \dd z^{\alpha} \right|
\end{align*}
We start by bounding the left hand-side from below.
\begin{align*}
\left|z^{\gamma} \sum_{\beta \in \tau} t^{-\nu(\beta)} \dd
z^{\beta-\gamma} \right| & = \left| \sum_{\beta \in \tau}
t^{-\nu(\beta)} z^\gamma \sum_{i=1}^{n} (\beta_{i} -\gamma_{i})
z^{\beta-\gamma-e_i} dz_{i} \right| \\ & = \sum_{i=1}^{n} \left|
\sum_{\beta \in \tau} t^{-\nu(\beta)} (\beta_{i} -\gamma_{i})
z^{\beta-e_i} dz_{i} \right| \\ & = \sum_{i=1}^{n} \left|
\sum_{\beta \in \tau} (\beta_{i} -\gamma_{i}) t^{-\nu(\beta)}
z^{\beta} \right| \\ & \geq  \frac{\left| t^{-\nu(\delta)}
z^{\delta}\right|}{\rho}  \end{align*} 
To see the last step, we interpret the penultimate line as the norm of a vector in $\bC^{n}$ which is the image of the $k$--vector $(t^{-\nu(\beta) } z^{\beta})|_{\beta \neq \gamma}$ under the linear map which takes the standard basis vectors to the vectors $(\beta - \gamma )_{\beta \in \tau}$. The constant $\rho$ was chosen to be a bound for the length distortion of this linear map.  It suffices therefore to bound the norm of $(t^{-\nu(\beta) } z^{\beta})|_{\beta \neq \gamma}$.  We have simply used the size of the $\delta$ component as the lower bound for the norm of this vector.

We now successively bound all the terms in the right hand-side.  We begin with the first term using \fullref{expo_decay}.
\begin{align*}
\left|\sum_{\alpha \in A} t^{-\nu(\alpha)} \phi_{\alpha}(z)
\dd z^{\alpha} \right| & \leq \sum_{\phi_{\alpha}(z) \neq 0}
t^{-\nu(\alpha)} \phi_{\alpha}(z) \left| \dd z^{\alpha} \right| 
\\ & \leq \sum_{\phi_{\alpha}(z) \neq 0} t^{-\nu(\alpha)} \sum_{i=1}^{n}
\left| \alpha_{i} \right| \left|z^{\alpha}\right| 
\\ & \leq N \left|A\right| e^{-c \epsilon \log(t)} \left| t^{-\nu(\delta)}
z^{\delta} \right| \\ & < \frac{ \left| t^{-\nu(\delta)}
z^{\delta} \right|}{5 \rho} \end{align*}
The same bound works for the last term.  We use \fullref{bound_f_t} to bound the third term.
\begin{align*}
 \left| \sum_{\beta \in \tau} t^{-\nu(\beta)} z^{\beta-\gamma}
\dd{z^{\gamma}}\right| & = \sum_{i=1}^{n} \left| \sum_{\beta \in
\tau} t^{-\nu(\beta)} z^{\beta-\gamma} \gamma_{i} z^{\gamma-e_{i}}
dz_{i} \right| \\ & = \sum_{i=1}^{n} |\gamma_i | \left| \sum_{\beta \in
\tau} t^{-\nu(\beta)} z^{\beta} \right| 
\\& < \frac{\left| t^{-\nu(\delta)} z^{\delta} \right|}{5 \rho}
\end{align*}
Finally, we bound the second term.
$$
\eqalignbot{
\left| \sum_{\alpha \in A} t^{-\nu(\alpha)} \dd \phi_\alpha(z)
z^\alpha \right| & \leq \sum_{\phi_{\alpha}(z) \neq 0} \left|
t^{-\nu(\alpha)} z^\alpha \right| \left| \dd \phi_\alpha(z) \right|\cr
& \leq \sum_{\phi_{\alpha}(z) \neq 0 } \left| t^{-\nu(\alpha)} z^\alpha
\right| \sum_{i=1}^{n} \left|\frac{\dd \phi_\alpha (z)}{\dd \log
\left| z_{i} \right|} \frac{dz^i}{\bar{z_{i}}} \right| \cr 
& = \sum_{\phi_{\alpha}(z) \neq 0} \left| t^{-\nu(\alpha)} z^\alpha \right|
\sum_{i=1}^{n} \left| \frac{\dd \phi_\alpha(z)}{\dd u_i} \right| \cr 
& \leq \sum_{\phi_{\alpha}(z) \neq 0 } \left|
t^{-\nu(\delta)} z^{\delta} \right| e^{-c\epsilon \log(t)}
\frac{4}{\epsilon \log(t)} \cr 
& \leq \frac{ \left|
t^{-\nu(\delta)} z^{\delta} \right| }{10 \rho}}\proved
$$
\end{proof}
Note that we can use the above bounds to conclude the following:
\begin{cor}\label{df-bound}
Under the hypotheses used in this section we have the bound $|\dd f_{t,s}| >
\frac{ \left| t^{-\nu(\delta)} z^{\delta} \right| }{10 \rho}$ .  \noproof
\end{cor}

\subsection{Moser's argument}

Our goal in this section is to prove the following:
\begin{prop} \label{pair}
 After rescaling the symplectic form, the pairs $(\ctorus{n},M_{t,s})$ are symplectomorphic for all $t$ and all $s$ for which they are defined.
\end{prop}

In the data of a pair $(N,M)$ of symplectic manifolds, we remember the
inclusion of $M$ as a symplectic submanifold of $N$.  In particular, a
symplectomorphism of pairs is a symplectomorphism between the total
spaces which carries the submanifolds to each other.  If we fix $t$,
we will define a map
\[ \psi_{s}\co \left( \ctorus{n},M_{t,0}) \to (\ctorus{n},M_{t,s} \right). \]
We will only prove the result for varying $s$, as the result for
varying $t$ is entirely analogous, although it requires rescaling the symplectic form by $\log(t)$.  We begin by showing that the
hypersurfaces $M_{t,s}$ are symplectomorphic using the fact that
our proof that $M_{t,s}$ is a symplectic hypersurface in fact proves
that $f_{t,s}$ is a symplectic fibration in a neighbourhood of each
$M_{t,s}$.  In particular we obtain a symplectic connection.  Let $Y_{t,s}$ be the
horizontal lift of $-\frac{\dd f_{t,s}}{\dd s}$.  We need to bound
$Y_{t,s}$ in order to show that it integrates to a flow.  But first,
we would like to obtain an explicit formula for it.  This is a problem
in linear algebra, for which we need a little notation.

Let $(V,\omega)$ be a symplectic vector space. The symplectic form $\omega$
determines an isomorphism between $V$ and $V^{*}$ which we denote by $\omega^{\#}$ and define by the formula
\[ \langle \omega^{\#}(v), w \rangle = \omega(v,w) .\]
The inverse isomorphism will be written $\omega_{\#}$.  We will also use the same symbols when we extend these maps $i$--linearly to the complexification of $V$ and its dual.  We equip $\bC$ with the standard symplectic form
\[ dx \wedge dy = \frac{i}{2} dz \wedge d\bar{z} .\]

\begin{lem}\label{splitting}
Let $f_{*}$ be a linear map between the symplectic vector spaces $(V, \omega)$
and $(\bC, \omega_{0})$ with symplectic kernel and with dual $f^{*}$. The
horizontal lift of the vector $a \in \bC$ is given by
\[f^{!}( a) =  \frac{2i}{\omega(f^* dz,f^{*}d\bar{z})} \omega_{\#}(f^{*}\omega_{0}^\# a) .\]
\end{lem}

\begin{proof}
We will work with the complexification of the vector space
throughout. First, note that the image of $\bC$ under the composition
$ \omega_{\#} \circ f^{*} \circ \omega_{0}^\#$ is precisely the
orthogonal complement of the kernel of $f_*$.  Hence it suffices to
show that the right hand side of our formula maps to $a$ under the
linear map $f_*$.  By linearity, it suffices to check this for
$\frac{\dd}{\dd z}$ and $\frac{\dd}{\dd \bar{z}}$, where, keeping in
mind our application, we've used a suggestive notation for the
standard basis of the complexification of $\bC$.  We check it for
$\frac{\dd}{\dd z}$ by computing the pairing with $dz$ and $d
\bar{z}$.
\begin{align*}
\bigg{\langle} f_{*} \left( \frac{2i}{\omega \left(f^*dz ,f^{*} d\bar{z} \right)}
\omega_{\#} f^{*} \left(\omega_{0}^\# \frac{\dd}{\dd z} \right) \right),dz \bigg{\rangle} & =
\frac{2i \big{\langle} f_{*}\left(
\omega_{\#}f^{*}\left( \frac{id \bar{z}}{2}\right)\right), dz \big{\rangle}  }{\omega(f^* dz,f^{*}d\bar{z})}
\\ & = \frac{- \big{\langle} \omega_{\#} f^{*}d \bar{z} , f^{*} dz \big{\rangle} }{\omega(f^* dz,f^{*}d\bar{z})}
\\ & = 1 
\\ \bigg{\langle} f_{*} \left( \frac{2i}{ \omega(f^* dz ,f^{*} d\bar{z} )} \omega_{\#}f^{*}\left( \omega_{0}^\# \frac{\dd}{\dd z} \right) \right),d\bar{z} \bigg{\rangle} & = \frac{ \big{\langle} \omega_{\#}f^{*}d \bar{z}, f^{*} d{\bar z} \big{\rangle} }{\omega(f^*dz,f^{*}dz)}
\\ & = 0.
\end{align*}
Of course, the same computation works for $\frac{\dd }{\dd \bar{z}}$.
\end{proof}

\begin{lem}
$Y_{t,s}$ is a bounded vector field, hence integrates to a local
diffeomorphism near the fibre $M_{t,s}$.
\end{lem}
\begin{proof}
By the previous lemma, we have the expression
\[ f_{t,s}^{!} \left(\frac{\dd f_{t,s}}{\dd s} \right) =  \frac{2i}{\omega(f_{t,s}^* dz,f_{t,s}^{*}d\bar{z})} \omega_{\#}f_{t,s}^{*}\left( \omega_{0}^\#\frac{\dd f_{t,s}}{\dd s} \right) .\]
Since $f_{t,s}$ is almost holomorphic, this vector
field has norm less than or equal to
\[ \frac{4 |\frac{\dd f_{t,s}}{\dd s}| \cdot |f_{t,s}^* dz|}{ |f_{t,s}^* dz|^{2}}. \]
By \fullref{df-bound}
\[ |f_{t,s}^* dz| \geq \frac{|\partial f_{t,s} |}{2} \geq  \frac{\left| t^{-\nu(\delta)} z^{\delta}\right|}{20 \rho } \]
with $\delta$ as in the proof of \fullref{family}.  On the other hand
\begin{align*}
\left| \frac{\dd f_{t,s}}{\dd s} \right| & = \left| \sum_{\phi_{\alpha}(z) \neq 0} \phi_{\alpha}(z) t^{-\nu(\alpha)} z^{\alpha} \right| 
\\ & \leq \frac{\left|  t^{-\nu(\delta)} z^{\delta}\right|}{5 \left| A \right| \rho N }
\end{align*}
as in \fullref{bound_f_t}.
\end{proof}

\begin{lem}
The flow $\bar{\psi}_{s}$ of $Y_{t,s}$ restricts to a
symplectomorphism between $M_{t,0}$ and $M_{t,s}$.
\end{lem}
\begin{proof}
First, we observe that by construction of $Y_{t,s}$
\begin{align*} 
\frac{\dd}{\dd s} f_{t,s} \left( \bar{\psi}_{s}(z) \right) & = {f_{t,s}}_{*}\left(
Y_{t,s}|_{\bar{\psi}_{s}(z)} \right) + \frac{\dd f_{t,s} }{\dd
s} \bigg{|}_{\bar{\psi}_{s}(z)} \\ & = 0.
\end{align*}
So the flow of $Y_{t,s}$ preserves the level sets of the function
$f_{t,s}$, and hence maps the zero fibres to each other.  To check
that this is, indeed, a symplectomorphism, it suffices to compute
\begin{align*}
\frac{\dd}{\dd s} \bar{\psi}_{s}^{*}\left( \omega_{M_{t,s}} \right) & =
\frac{\dd}{\dd s} \bar{\psi}_{s}^{*} \left( \omega|_{M_{t,s}} \right)\\ & =
\bar{\psi}_{s}^{*}\left( \cL_{Y_{t,s}} \omega|_{M_{t,s}}\right) \\ & =
\bar{\psi}_{s}^{*}\left(d(i_{Y_{t,s}} \omega)|_{M_{t,s}}\right) \\ & = 0
\end{align*}
using Cartan's formula, the fact that $\omega$ is closed, and that
$Y_{t,s}$ is orthogonal to $M_{t,s}$.  In particular, $i_{Y_{t,s}}
\omega$ vanishes on $M_{t,s}$ and hence also its exterior derivative.
\end{proof}

We must now check whether this symplectomorphism can be embedded in an
appropriate symplectomorphism of pairs.  Similar results are well
known in the compact case \cite{ruan} and the proof extends to the
open case provided one can bound the vector fields.

\begin{lem}\label{flow_integrates}
There exists a Hamiltonian time-dependent vector field $Y'_{t,s}$ on
$\ctorus{n}$ which is supported in a neighbourhood of $M_{t,s}$ and
which integrates to a symplectic flow $\psi_{s}$ that maps $M_{t,0}$
to $M_{t,s}$.
\end{lem}
\begin{proof}
There exists a function $H_{s}$ which vanishes on $M_{t,s}$, and such
that $dH_{s}|_{M_{t,s}} = i_{Y_{t,s}} \omega$.  This can be seen
locally by taking any trivialization of a neighbourhood of $M_{t,s}$
and thinking of $ i_{Y_{t,s}} \omega$ as an initial value problem for
a function, then using a partition of unity to produce a global
function.  We can always choose $H_{s}$ to be supported in a
neighbourhood of $M_{t,s}$.  The form $dH_{s}$ gives the closed form
which corresponds to our vector field $Y'_{t,s}$.  Since our initial
values are bounded, an appropriate choice of cutoff function ensures
that the vector field $Y'_{t,s}$ is also bounded, hence integrates to
a Hamiltonian isotopy.
\end{proof}
\begin{rem}
Note that we can choose the germs of the vector fields $Y'_{t,s}$ and $Y_{t,s}$ to be $C^{\infty}$ close at $M_{t,s}$.  In particular, this means we can assume that $\psi_s$ respects the fibration to arbitrarily high order at $M_{t,s}$. 
\end{rem}

Note that if $J$ is any compatible complex structure $J$ on $V$ such
that $f_{*}$ is a map of complex vector spaces, then the symplectic
horizontal lift agrees with the orthogonal lift which for a unit
length vector $v$ is given by
\[ f^{!}(v) =  \frac{1}{g\left(f^{*}g_{0}^{\#}v,f^{*}g_{0}^{\#}v\right)} g_{\#}f^{*}g_{0}^\# v.\]
In particular, if $f_*$ is the derivative at the origin of a
symplectic fibration $f\co N \to \bC$, we consider the real valued
function $f_{\bR}$, which is the restriction of $f$ to the inverse
image of the real axis. Since $df_{\bR} = f^{*}g_{0}^{\#}\frac{\partial}{\partial x}$, we conclude
\[ \grad{f_{\bR}} = |df_{\bR}|^2 f^{!}\frac{\partial}{\partial x}.\]
In future section, we will let $\Lambda$ be a compact Lagrangian submanifold of the zero fibre of $f$, which we denote $M$, and consider parallel transports of $\Lambda$ along straight half-rays originating at the origin, and moving in the left half-plane.  We write such a line as 
\[ \gamma_{c}(t) = - t (1+2ci). \]
Denote the parallel transport along this line by $L_{c}$, and the distance to $\Lambda$ by $r_{\Lambda}$.  On $L_{0}$, the gradient of $r_{\Lambda}$ and the lift of $\frac{\partial}{\partial x}$ are both orthogonal to $\Lambda$, and satisfy
\[ \grad r_{\Lambda} =   |df_{\bR}| f^{!}\frac{\partial}{\partial x}.\]
\begin{lem} \label{sec}
Let $J$ be an almost complex structure in a neighbourhood of $M$ for
which $f$ is holomorphic and $M$ is an almost complex submanifold.  If
we identify a neighbourhood of $L_0$ with its cotangent bundle as in
\fullref{weinst}, then $L_c$ is given as the graph of a $1$--form
which agrees with the differential of $c r^{2}_{\Lambda}(u)$ to first order in $u$.
\end{lem}
\begin{proof}
Let $L'_{c}$ be the graph of the differential of $c r^{2}_{\Lambda}$. Since the differential of $cr^{2}_{\Lambda}(u)$ vanishes on $\Lambda$, $L_c$ and $L'_c$ have the
same boundary.  To show that the $1$--forms that define them agree to
first order in $u$, it suffices to show that
their tangent spaces agree on the boundary.  However, $TN_{p}$
decomposes as an orthogonal direct sum $TM \oplus \bC$, and it is clear
that
\[ TL_{c} = T \Lambda \oplus \bR \cdot f^{!}\left( \frac{\partial}{\partial x} + 2c \frac{\partial}{\partial y} \right) ;\]
 so it remains to compute $TL'_{c}$.  To do this, we observe that
 $r_{\Lambda}$ vanishes on $\Lambda$, and that since $TM$ is closed
 under $J$, and $J$ is compatible with $\omega$, its orthogonal
 complement is also closed under $J$.  We may therefore ignore the
 directions corresponding to $T \Lambda$ in our computation.  Now the
 orthogonal complement of $T \Lambda$ in $TL_{0}$ is spanned by $f^{!}
 \frac{\partial}{\partial x}$, so the component of the tangent space
 of $TL'_c$ which lies in the orthogonal complement to $T \Lambda$ is spanned
 by
 \begin{equation}\label{normal} f^{!} \frac{\partial}{\partial x} + \nabla_{f^{!} \frac{\partial}{\partial x}} \left(2cr_{\Lambda}(u) J   \grad r_{\Lambda} (u) \right). \end{equation}
Consider the vector field
\[ Z =  J\left( \grad r_{\Lambda} -  |df_{\bR}| f^{!}\frac{\partial}{\partial x} \right) .\]
We know that this vector field vanishes at $\Lambda$.  In addition, we know that the vector field $J f^{!}\frac{\partial}{\partial x} = f^{!}\frac{\partial}{\partial y}$ is the lift of a geodesic flow in the base.  We conclude that on $\Lambda$ the second term can be computed as follows:
\begin{align*}
 \nabla_{f^{!} \frac{\partial}{\partial x}} \left( r_{\Lambda}(u) J   \grad r_{\Lambda} (u) \right) & = \nabla_{f^{!} \frac{\partial}{\partial x}} \left( r_{\Lambda}(u)  |df_{\bR}| J f^{!} \frac{\partial}{\partial x}  + r_{\Lambda}(u) Z \right)
\\ & =   \nabla_{f^{!} \frac{\partial}{\partial x}} \left( r_{\Lambda}(u)  |df_{\bR}| f^{!} \frac{\partial}{\partial y} \right) + r_{\Lambda}(u) \nabla_{f^{!} \frac{\partial}{\partial x}} Z 
\\ & \quad +  f^{!} \frac{\partial}{\partial x}(r_{\Lambda}(u)) Z
\end{align*}
Since $r_{\Lambda}$ and $Z$ vanish on $\Lambda$, it remains to compute that:
\begin{align*}
 \nabla_{f^{!} \frac{\partial}{\partial x}} \left( r_{\Lambda}(u)  |df_{\bR}| f^{!} \frac{\partial}{\partial y} \right) & = f^{!} \frac{\partial}{\partial x} \left( |df_{\bR}| r_{\Lambda}(u) \right) \frac{\partial}{\partial y} 
\\ & = \left( r_{\Lambda}(u) f^{!} \frac{\partial}{\partial x} (  |df_{\bR}|) + |df_{\bR}| f^{!}  \frac{\partial}{\partial x}( r_{\Lambda}(u))  \right)  f^{!} \frac{\partial}{\partial y}
\\ & =  \grad r_{\Lambda}(u)( r_{\Lambda}(u))   f^{!} \frac{\partial}{\partial y}
\\ & =   f^{!} \frac{\partial}{\partial y}
\end{align*} In the last line, we've used the fact that the restriction
gradient of $r_{\Lambda}$ to $\Lambda$ has unit norm.  Plugging this back in Equation \eqref{normal} completes the proof of the lemma.
\end{proof}
\subsection{The zero-section} \label{zero}
In this section, we construct an admissible Lagrangian submanifold of
$\ctorus{n}$ with boundary on $M_{t,1}$, which will be the first step of future constructions. The construction can be done for any Laurent polynomial $f$.  However, we will specialize to the situation where we're considering the mirror $W$ of a smooth toric variety.  In particular, 
 \[W_{t,1}=-1+\sum_{0 \neq \alpha \in A} t^{-\nu(\alpha)} (1 - \phi_{\alpha} (z) )z^{\alpha}.\] 
Consider the $\epsilon \log(t)$ neighbourhood of the $n-2$ skeleton of the polytope $Q_t = C_{0,t}$.  The intersection of the amoeba of $W_{t,1}$ with a neighbourhood of $Q_t$ agrees with the boundary of $Q_t$ outside this neighbourhood of the $n-2$ skeleton.

Recall that by thinking of $\ctorus{n}$ as the cotangent bundle of
$\mathbb R^n$ modulo the lattice $\mathbb{Z}^{n}$, we identified a
natural zero-section.  One might equivalently think of the inclusion
of $(\mathbb R^{+})^{n}$ into $\ctorus{n}$.

\begin{lem} \label{boundary_zero}
Near the component of the complement of the amoeba corresponding to
the origin, the zero section intersects $M_{t,1}$ in a smooth Lagrangian
sphere of $M_{t,1}$. This Lagrangian sphere is the boundary of a compact
subset of the zero section which is diffeomorphic to the $n$--ball.
\end{lem}
\begin{rem}
It is easy to see that the zero section does intersect the zero level
set of the linear function
\[  -1 + \sum_{i=1}^{n}  t^{-\nu(i)} z_{i} \]
in a smooth Lagrangian submanifold.  Hence, by the multiplicative
change of coordinates, the same holds for
\[g(z) = -1 + \sum_{0 \neq \alpha \in \tau} t^{-\nu(\alpha)} z^{\alpha}\]
whenever $\tau$ is a minimal simplex.
\end{rem}
\begin{proof}
We use the polyhedral decomposition provided by the tropical limit, and prove the result for $M_{t,1}$.  Consider a point $z$ such that $\Log(z) \in O_{\check{\tau}}$. If $\tau$ is a $k$--simplex, we know that all but $k+1$ terms in
 \[W_{t,1}= -1 + \sum_{0 \neq \alpha \in A} t^{-\nu(\alpha)} (1 - \phi_{\alpha} (z) )z^{\alpha} \] 
are zero. It suffices therefore to check the result for a simplex as in the previous remark.  So consider
\[g_{k,s}(z) = -1 + \sum_{0 \neq \alpha \in \tau} t^{-\nu(\alpha)}(1 - s \phi_{\alpha} (z) ) z^{\alpha}\]
which we think of as a $1$--parameter family of functions defined on a
domain where $g_{k,1}$ agrees with $W_{t,1}$.  By the previous remark,
it will be sufficient to show that reality is preserved by our flow,
which reduces to showing that the vector field which defines the flow
is tangent to $(\mathbb R^{+})^{n}$ on the intersection of $(\mathbb
R^{+})^{n}$ with the zero level set.  Using the result of \fullref{splitting}, and the expression of the symplectic form on
$\ctorus{n}$ in standard coordinates, we are reduced to showing that
\[ g_{k,s}^{*}(\omega_{0}^\# \frac{\partial g_{k,s}}{\partial s}) = \sum_{i=1}^{n} h_{i}(z) dy_i \]for some real functions $h_i$.  Since $g_{k,s}$ restricts to a real function on the intersection of  $(\mathbb R^{+})^{n}$ with the zero level set, the restriction of $\frac{\partial g_{k,s}}{\partial s}$ to this intersection is some real multiple of $\frac{\partial}{\partial x}$, hence $\omega_{0}^\# \frac{\partial g_{k,s}}{\partial s}$ is some real multiple of $dy$.  The result would therefore follow from showing that for all $1 \leq i \leq n$, the following derivative of the imaginary part is zero:
\[ \frac{\partial \im( g_{k,s})}{\partial x_{i}} = 0 \]
But in fact, the function $\im( g_{k,s})$ itself vanishes on the real
locus, so this establishes that the flow maps the positive real locus
of $g_{k,0}$ into the positive real locus of $g_{k,1}$.  Since we can
also use the reverse flow, this map gives a diffeomorphism between the
two sets.

To see that the result is simply an $n$--ball with boundary, we note
that the boundary of the component of the complement of the amoeba
corresponding to the origin is a polytope which is itself
homeomorphic to a sphere which bounds $\cL_{\infty}$.  It is clear from the local model that $\cL_{\infty}$ is star-shaped about the origin, so it must be diffeomorphic to the standard ball.
\end{proof}

We will denote this Lagrangian ball with boundary on $M_{t,1}$ by $
\log(t) \cL$. 

\begin{rem} \label{switch}
Note that the results of the previous section allow us to conclude that any
construction which is performed for admissible Lagrangians with
respect to $W_{t,1}$ can be ``pulled back" to the original complex
hypersurface $M$.  Hereafter, we will change our notation and
write $M$ for the image of $M_{t,1}$ under the diffeomorphism of
$\ctorus{n}$ induced by the diffeomorphism of the base
\[ u \to \frac{u}{\log(t)} .\]
 This induces the rescaling of the symplectic form which appears in the statement of \fullref{pair}.  But rescaling preserves symplectic and Lagrangian
 submanifolds, so our constructions are insensitive to it.  In particular, we produce a Lagrangian ball $\cL$ with boundary on $M$.

We will also use $W$ for the pull back of $W_{t,1}$ under the above ``conformal'' symplectomorphism. We will still perform most computations using the coordinates $z_i$ on $\ctorus{n}$ as a subset of $\bC^{n}$.  Strictly speaking, these differ by a ``rescaling'' from the coordinates used in the previous section.  These coordinates come with a natural complex structure which will use.
\end{rem}

There are two reasons that allow us to use this integrable almost complex structure.  The first reason is that we proved in \fullref{cup_product} that the holomorphic triangles that we will be considering are already regular.  This means that there is no issue in the interior.  More problematic is the fact that we will also be considering the distance function with respect to the usual metric in order to approximate parallel transport as in \fullref{sec}.  However, it is easy to see from the remark following \fullref{flow_integrates} that we can choose the usual complex structure to be $C^1$--close near $M$ to an admissible almost complex structure $J$. This implies that the distance function with respect to the usual metric yields a $C^1$--close approximation to parallel transport with respect to $J$, which is all we will need.

\section{The admissible Lagrangians $\cL(j)$} \label{admiss-sect}
\subsection{Desiderata} \label{des}
Consider the symplectic hypersurface $M$ with amoeba $\cA$.  By
\fullref{close-amoeba}, $\cA$ is $\epsilon$--close to the polyhedral complex $
\Pi$.  If we focus only on $Q$, we know that the amoeba agrees with
the polyhedral complex away from an $\epsilon$ neighbourhood of its
$n-2$ skeleton.  In the next two sections, we will construct a
Hamiltonian function $H$ which should be thought of as a modification
of $H_{\infty}$ and whose time--$1$ flow preserves $M$.  But first, we
would like to describe the properties that $H$ must satisfy, and show
that all possible choices of $H$, up to admissible Hamiltonian isotopy, are
equivalent.  

Continuing with the change in notation introduced in the last section, we have a cover $O_{\sigma}$ of $\bR^{n}$ indexed by the faces $\sigma$ of $\Pi$ such that each $O_{\sigma}$ is a neighbourhood of a large open subset of $\sigma$.  Recall that we have been using the
Euclidean metric on $\bR^{n}$ to identify each cotangent fibre with a
copy of $\bR^{n}$; let $\pi$ denote this cotangent fibration.

\begin{defin}
A Lagrangian submanifold $\Lambda$ of $M$ is a {\em boundary for
$\cL(1)$} if $\Lambda$ lifts to a submanifold $\tilde{\Lambda} \subset
T \bR^{n}$ such the following properties are satisfied:
\begin{itemize}
\item $\tilde{\Lambda}$ is a section of $T \bR^{n}|_{\partial
\cL}$, and
\item if $z \in \tilde{\Lambda}$ and $\pi(z) \in \partial \cL \cap O_{\sigma}$ then $z$ lies in the affine subspace supporting $-2\pi \sigma$, and
\item $\tilde{\Lambda}$ is the restriction to $\cL$ of an exact section of
$T^{*} \bR^{n}$.
\end{itemize}
\end{defin}

Note that the last condition is vacuous if $n>2$ since $H^{1}(\partial
\cL)$ vanishes in this case.

\begin{lem} \label{boundary_exists}
Given $\Lambda$, a boundary for $\cL(1)$, and a vector $v$ in the left
half plane, there exists an admissible Lagrangian $\cL(1)$ which is an
exact section of $T \bR^{n}|_{\cL}$ with boundary equal to
$\Lambda$ and whose image under $W$ agrees near the origin with a curve $\gamma$ whose tangent vector at the origin is $v$.  Furthermore, any two such
Lagrangians, for (possibly) different tangent vectors at the origin,
are isotopic through admissible Lagrangians.
\end{lem}
\begin{proof} To prove the first statement, we prove that the parallel transport of $\Lambda$ along a short curve $\gamma$ with tangent vector $v$ at the origin yields a section of $T \bR^{n}|_{U}$ where $U$ is a $1$--sided neighbourhood of $\partial \cL$ in $\cL$.  As usual, we consider the problem on each element of the cover $O_{\sigma}$.   Note that the condition that $\tilde{\Lambda}$ lies in the hyperplane that supports $-2\pi \sigma$ implies that $\tilde{\Lambda}$ is locally the product of a section of $T \sigma$ with the zero section of a transverse symplectic slice.  By \fullref{semi-local} we know that in $O_{\sigma}$ the function $W$ is independent of the directions corresponding to $\sigma$, so parallel transport only depends on the directions which lie in the transverse slice, which we may locally identify with a lower dimensional torus. 
\begin{claim}
If $\cL_v$ is the parallel transport of the zero-section along a sufficiently short curve $\gamma$ starting at $0$ with tangent vector $v$ in the left half-plane, then $\pi$ projects $\cL_v$ homeomorphically onto a neighbourhood of the boundary of $\cL$.
\end{claim}
\begin{proof}[Proof of Claim]
In the proof of \fullref{boundary_zero} we noticed that $W_{t,1}$ takes negative real values on $\cL$. In particular, this implies that $\cL$ is the parallel transport of its boundary along the negative real axis.  The result now follows immediately from the proof of \fullref{sec}.
\end{proof}

Constructing a Lagrangian $\cL_v(1)$ as the parallel transport of $\Lambda$ along such a curve $\gamma$ yields locally a section $T \bR^{n}|_{\cL}$.  Our third condition on $\Lambda$ guarantees that this section is exact, hence any extension of the function that defines it will yield the desired section of $T^{*} \bR^{n}|_{\cL}$.

To prove uniqueness, we observe that any Lagrangian $\cL(1)$ is the graph
of the differential of a function $H$ with prescribed derivatives at
the boundary. The linear interpolation between two such functions $H$
and $H'$ defines an isotopy which may not preserve admissibility. But
choosing a sufficiently fine cover of the interval $[0,1]$ as in \fullref{invariance} allows us to modify this isotopy to one in which
admissibility is preserved.
\end{proof}

This means that we can unambiguously define $\cL(1)_{\Lambda}$ to be
an admissible exact section of $T^{*} \bR^{n}|_{\cL}$ whose boundary
is $\Lambda$.  We would like to eliminate the dependence on $\Lambda$.
Recall from \fullref{semi-local} that $W$ is independent of the
directions corresponding to the tangent space of $\sigma$.  In
particular, the proof of the previous lemma together with the
convexity of the faces of $\sigma$ imply that we can use interpolation as in \fullref{invariance} to prove the following:

\begin{lem}
If $\Lambda$ and $\Lambda'$ are boundaries for $\cO(1)$, then
$\cL(1)_{\Lambda}$ and $\cL(1)_{\Lambda'}$ are isotopic through
admissible Lagrangians. \noproof
\end{lem} 

So we can define admissible Lagrangian submanifolds $\cL(1)$ which are
well defined up to Hamiltonian isotopy.  In the next sections, we will
produce an explicit model for $\cL(1)$ that will allow us to compute
its Floer cohomology groups.  We also note that this construction
naturally extends to the construction of Lagrangians $\cL(j)$, for any
integer $j$, whose lifts lie in $-2 \pi j \sigma$ in every set $O_{\sigma}$.

\subsection{A boundary for $\cL(1)$}

In this section, we will construct a boundary for $\cL(1)$.   First, we must choose an appropriate cover of a neighbourhood of $\Pi$.  Given $\sigma$ a facet of $\Pi$, Let $U_{\sigma}$ be its $\epsilon$--neighbourhood and $V_{\sigma}$ its $2 \epsilon$ neighbourhood. If we then consider the $n-2$ skeleton, there exists a constant, $c_{n-2}$ such that every point in the $2 \epsilon$ neighbourhood of $\Pi$, but which is distance more than $c_{n-2} \epsilon$ from the $n-2$ skeleton is closest to a unique facet of $\Pi$.   The value of $c_{n-2}$ is independent of $\epsilon$, and is related to the ``angles'' between the facets near their intersections.  Given a face $\sigma$ in the $n-2$ skeleton, let $U_{\sigma}$ be the $c_{n-2} \epsilon$ neighbourhood of $\sigma$, and $V_{\sigma}$ the $2 c_{n-2} \epsilon$ neighbourhood of $\sigma$. 

Repeating this process inductively, we obtain constants $c_{i} \geq 1$ for each $i \leq n-1$ independent of $\epsilon$, and define $U_{\sigma}$ to be the $c_i \epsilon$ neighbourhood of $\sigma$, and $V_{\sigma}$ its $2 c_i \epsilon$ neighbourhood for $i = \dim(\sigma)$.  Let 
\begin{align*}
U_{i} & = \bigcup_{\dim(\sigma) = i} U_{\sigma}
\\ V_{i} & = \bigcup_{\dim(\sigma) = i} V_{\sigma}.
\end{align*}
The main property satisfied by this cover is the following:
\begin{lem}
Every point of $V_{i}-U_{i-1}$ is closest to a unique face of dimension $i$. \noproof
\end{lem}

\begin{figure}
\begin{center}
\includegraphics[width=1.5in]{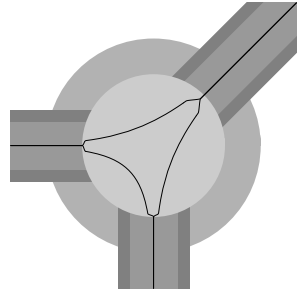}
\end{center}
\caption{The sets $U_0$, $V_0$, $U_1$ and $V_1$ for the hyperplane in
$\ctorus{2}$ shown in progressively darker shades of gray}
\label{amoebas}
\end{figure}

In the interior of each $i$--face $\sigma$, we identify the normal bundle of $\sigma$ with a neighbourhood, and define a self map $\psi_{\sigma}$ of this bundle which is radial, collapses $U_{\sigma}$ to $\sigma$, and is the identity away from $V_{\sigma}$.  We would like to take the composition of all $\psi_{\sigma}$ over the $n-2$ skeleton, to obtain a self map of $\bR^{n}$.  There are two issues with this:
\begin{itemize}
\item $\psi_{\sigma}$ does not extend continuously to $\bR^{n}$
because $\sigma$ has non-empty boundary (if $\dim(\sigma) \neq 0$).
\item The supports of $\psi_{\sigma}$ for different faces are not
disjoint.
\end{itemize}

We resolve both of these issues by defining $\psi_{i}$ to be the
composition of $\psi_{\sigma}$ for $\sigma$ a face of dimension $i$ in
an arbitrary order, and defining $\psi$ to be the composition
\[ \psi = \psi_{0} \circ \psi_{1} \circ \cdots \circ \psi_{n-1} .\]
The key point is that the set where $\psi_{i}$ is not continuous, or
depends on the order chosen, lies entirely within $U_{i-1}$, which is
collapsed by the next map ($\psi_{i-1}$) to the $i-1$ skeleton.  It
should therefore be clear from the construction that $\psi$ is a well
defined smooth self map of $\bR^{n}$ which does not depend on the
order of compositions chosen for $\psi_{\sigma}$.  Further, we have
the following:
\begin{lem}
The map $\psi$ can be chosen to satisfy the following additional properties:
\begin{itemize}
\item $\psi$ preserves the $i$--skeleton of $\Pi$.  In particular, it
maps every face to itself.
\item $\psi$ maps every point $p \in U_{i}$ to the nearest $i$--face of $\Pi$.  If $p \in U_{i} - V_{i-1}$ then $\psi(p)$ is the nearest point on the nearest 
$i$--face.
\item $\psi$ maps the amoeba $\cA$ to the tropical amoeba $\Pi$.
\item $\psi$ can be chosen to satisfy the condition that its derivative is bounded independently of $\epsilon$.  Specifically, given any tangent vector $v$ at $u$, \[|\psi_{*}(u)(v)| < 4^{n} |v|.\]
\end{itemize}
\end{lem}
\begin{proof}
The first two properties are immediate from the construction, so we only prove the other two.  Note that we have $O_{\sigma} \subset U_{\sigma}$.  In particular, since the localization of the tropical amoeba lies entirely in $\cup_{\sigma} O_{\sigma}$, $\psi$ maps the localization of the tropical amoeba to the tropical amoeba.  Finally, since each $\psi_{\sigma}$ collapses a neighbourhood of radius $c_i \epsilon$, and is the identity outside a neighbourhood of radius $2 c_i \epsilon$, its differential can be chosen to expand lengths of vectors by a factor less than $4$.  The result follows by observing that in a neighbourhood of any point $p$, $\psi$ can be written as the composition of at most $n$ maps corresponding to the faces of each dimension which are closest to $p$.
\end{proof}

Consider the Hamiltonian function given by
\begin{equation*} \label{define_boundary} H^{\partial}(u) = -2 \pi {\Big{\langle}} u - \frac{\psi(u)}{2}, \psi(u) {\Big{\rangle}} - \pi a \langle \psi(u) , \psi(u) \rangle \end{equation*}
on the set $\cL \cap \cup_{i} U_i$ , and define $\Lambda$ to be the graph of $dH^{\partial}|_{\partial \cL}$.  To keep the notation simpler, we do not make the dependence on $a$ explicit.

\begin{lem}
$\Lambda$ is a boundary for $\cL(1)$.
\end{lem}
\begin{proof}
Let us first set $a=0$.  We compute that for any tangent vector $v$,
\begin{align*}
dH^{\partial}(u)(v) & = -2 \pi \Big{\langle} v - \frac{ \psi_{*}(u) (v)}{2}, \psi(u) \Big{\rangle}  -2 \pi \Big{\langle} u - \frac{\psi(u)}{2}, \psi_{*}(u)(v) \Big{\rangle}
\\ & =  \langle v , - 2\pi \psi(u) \rangle  - 2 \pi \langle u - \psi(u), \psi_{*}(u)(v) \rangle.
\end{align*}
The only condition we need to check is that if $u \in O_{\sigma}$, this co-vector lies in the hyperplane supporting $-2\pi \sigma$ after passing to the tangent space.  By the previous lemma, $\psi(u)$ lies in $\sigma$, so this is already true of the first term.  We must prove that the second term lies in the tangent space of $\sigma$.  In particular, it suffices to show that the corresponding $1$--form annihilates every tangent vector normal to $\sigma$.  But  in $U_{\sigma}$, $\psi_{\sigma}$ collapses the directions orthogonal to $\sigma$ so this is indeed true. 

To address the case $a \neq 0$, we observe that the term depending on $a$ has differential
\[ v \mapsto - 2 \pi a \langle  \psi_{*}(u) (v) , \psi(u) \rangle \]
which, again, corresponds to a term in the tangent space of $\sigma$.
 \end{proof}

\subsection{Parallel transport of $\Lambda$}
As in the proof of \fullref{boundary_exists}, we begin the construction of $\cL(1)$ by parallel transport of the boundary along a curve $\gamma \subset \bC$ starting at the origin and moving in the left half plane.  We proved that this parallel transport can be expressed locally as an exact section over $\cL$.  Recall that in order to compute
\[ HF^{*}(\cL, \cL(1)) \]
we need to ensure that $W(\cL(1))$ agrees with a curve $\gamma$ near the origin, and that the angle between the negative $x$--axis and $\gamma$ is between $0$ and $\frac{\pi}{2}$.  In particular, we would be considering $c>0$ in the notation of \fullref{sec}.  Let us write $\cL_{c}(1)$ for a Lagrangian that satisfies this condition of being a parallel transport of $\Lambda$ along a line of slope $2c$ at the origin.

\begin{lem}
The graph of $dH^{\partial}$ agrees with $\cL_{0}(1)$ near $\Lambda$. 
\end{lem}
\begin{proof}
We only need to prove that the image under $W$ of the graph of $dH^{\partial}$ near $\partial L$ is the negative real axis.   By construction, $dH^{\partial}$ lies in the hyperplane supporting $-2 \pi \sigma$, so we may think of it as the product of the zero section of $(\bC^{\star})^{\check{\sigma}}$ with a non-trivial section of $T \sigma$.  Since in $O_{\sigma}$, $W$ only depends on the directions transverse to $\sigma$, it suffices to prove that the zero-section $\cL$ is the parallel transport of its restriction to the boundary along the negative real axis.  We already observed this in the proof of \fullref{boundary_exists}.
\end{proof}

The graph of the  differential of
\[ H^{\partial}_{c} (u) = -2 \pi \left( \Big{\langle} u - \frac{\psi(u)}{2}, \psi(u) \Big{\rangle}+  \frac{a}{2} \langle \psi(u) , \psi(u) \rangle \right)+ c r_{\Lambda}^2(u) \]
 is therefore $C^1$--close to an admissible Lagrangian by \fullref{sec}.   In particular, using \fullref{C1-close}, we can construct an admissible Lagrangian $\cL_{c}(1)$ that is $C^{1}$--close to this graph near $\Lambda$, and that agrees with it away from an arbitrarily small neighbourhood of $\Lambda$.  In practice, this means that we can think of the graph of $dH^{\partial}_{c}$ as being admissible.

We will interpolate between $\cL_{c}(1)$ and the linear Lagrangians $\cL_{\infty}(1)$ constructed in \fullref{twist1}.  We will only do the harder case where $c > 0$ which corresponds to choosing a model of $\cL(1)$ such that $\cL$ and $\cL(1)$ are a positively oriented pair.

Let $\delta$ be an arbitrarily small positive number.  We assume that the origin is contained in the interior of $Q$.  Choose $\epsilon' >0$ such that $(1+\epsilon') Q$ is contained in a $2\delta$ neighbourhood of $Q$.  Note that $(1+\epsilon') \partial Q$ lies entirely outside of $Q$.  Further, given any face $\sigma$ of $\partial Q$, $(1+\epsilon') \sigma$ lies in the interior of $2 \sigma - Q$.  Choose $\epsilon >0$ and go through the constructions of the previous sections so that $(1+\epsilon')U_{\sigma}$ also lies in the interior of $2 \sigma - Q$.  Note that we can choose $\epsilon < \epsilon' < \delta$.

There exists a number $0< b < \frac{1}{2}$ such that $r_{\Lambda}$ has no critical points in the $b \epsilon$ neighbourhood of $\Lambda$.  Let $g$ be a non-negative cutoff function satisfying the following conditions:
\begin{align*}
g(u) = 0 & \Leftrightarrow  r_{\Lambda}(u) \leq b \epsilon
\\ g(u) = 1 & \Leftrightarrow  r_{\Lambda}(u) \geq \epsilon \end{align*}
\[ \left| dg (u) \right| \leq \frac{2}{\epsilon} \]
Consider the Hamiltonian $H_{c}$ given by the formula
\begin{align*}
& - (1-g(u))  \left( \pi \langle 2u - \psi(u), \psi(u) \rangle + (\epsilon' + 5 \epsilon \epsilon') \pi \langle \psi(u), \psi(u) \rangle - c r^{2}_{\Lambda}(u) \right) \\ & \qquad - g(u) \pi (1+ \epsilon') \langle u , u \rangle .
\end{align*}
Note that for $r_{\Lambda}(u) < b \epsilon$, this function is just $H^{\partial}_{c}$ for $a = \epsilon' + 5 \epsilon \epsilon'$.  In particular, the corresponding Lagrangian is $C^{1}$--close to $\cL_c(1)$.  Since the images of $\cL_c(1)$ and $\cL$ under $W$ do not intersect, it is clear that the Lagrangian corresponding to $H_c$ does not intersect $\cL$ in this region. Let us abuse notation and write $\cL_c(1)$ for the graph of the function $H_c$.   When $r_{\Lambda}(u) \geq \epsilon$, the Hamiltonian is
\[ - \pi (1+ \epsilon') \langle u , u \rangle \]
and the intersection points between $\cL$ and $\cL(1)_c$ in this region are exactly the $ \frac{1}{1+\epsilon'}$ integral points of $Q$.  
\begin{lem}
$dH_{c}$ takes no $2 \pi$--integral values in the region where $b \epsilon \leq r_{\Lambda}(u) \leq \epsilon$.
\end{lem}
\begin{proof}
Since it suffices to prove non-integrality of some coefficient, we use the cover $O_{\sigma}$ and consider only the factor transverse to $\sigma$.  In other words we may assume, after possibly reducing to a factor, that $u$ lies in the $\epsilon$ neighbourhood of a vertex $u_0$ of $Q$.  Note that this means that $\psi(u) = u_0$, and that the differential of $\psi$ vanishes.  We compute that in this case, $dH_{c}$ equals
\begin{align*}
& - 2 \pi u_0 + g(u) \left( 2 \pi u_0 - 2\pi (1+ \epsilon') u \right) + (1 - g(u)) 2 c r_{\Lambda}(u) dr_{\Lambda} (u)
\\& \quad + \left( \pi \langle 2u - u_0, u_0 \rangle + (\epsilon' + 5 \epsilon \epsilon') \pi \langle u_0 , u_0 \rangle - c r^{2}_{\Lambda}(u) - (1+ \epsilon') \pi \langle u , u \rangle \right) dg(u) .
\end{align*}
Since $-2 \pi u_0$ is an integral vector it suffices to show that the remaining terms do not add up to $0$, and have norm much less that $1$. The main observation is the following:
\begin{claim}
The inner product of any of the vectors $dg$, $dr_{\Lambda}$ or $u_0 - (1+ \epsilon')u$ with $ u_0$ is non-positive.
\end{claim}
\begin{proof}
For the last vector, this follows from the choices made about $\epsilon'$ and $\epsilon$ at the beginning of this section.  In particular, we can choose $\epsilon$ much smaller than $\epsilon'$, so this vector is well approximated by $ - \epsilon' u_0$.  Since $dg$ can be chosen approximately parallel to $dr_{\Lambda}$, we will simply prove the claim for this last vector.  But the gradient of $dr_{\Lambda}$ is a scalar multiple of the normal vector at the nearest point in $\Lambda$. On the other hand, it is easy to compute that the negative of the gradient of the Laurent polynomial
\[W_{t}(z) = -1 + \sum_{t} t^{-\nu(\alpha)} z^{\alpha} \]
satisfies this condition on $\Lambda$ (this is simply the statement that $W$ decreases as we move towards the interior of $\cL$).  Since the cutoff functions were chosen $C^{1}$--small with respect to $W_t$, this result is unaffected by the presence of $\phi_{\alpha}$.  We can then rescale by $\log(t)$ to work with our new conventions.
\end{proof}
Each of these vectors appears in the expression of $dH_c$.  Since the terms that they correspond to cannot all vanish at the same time, it suffices to show that the coefficients of each of these term is non-negative and bounded.  It is clear that the first two coefficients are indeed non-negative, and we have the bounds
\begin{align*}
\left|  g(u) ( 2 \pi u_0 - 2 \pi (1+ \epsilon') u) \right| & \leq 2 \pi \delta
\\ \big{|} (1 - g(u)) 2 c r_{\Lambda}(u) dr_{\Lambda} (u) \big{|} & \leq 2 c \epsilon.
\end{align*}
In the second bound, we used the fact that $|d r_{\Lambda} | \leq 1$.  It remains, therefore, to consider the coefficient of the last term. We re-write it as
\[ \epsilon'  \pi \big{(} (1+ 5 \epsilon) \langle u_0 , u_0 \rangle - \langle u, u \rangle \big{)} - \pi \big{(} \langle u , u \rangle - 2\langle u , u_0 \rangle + \langle u_0 , u_0 \rangle \big{)} - c r_{\Lambda}^2(u),\]
which simplifies to
\[ 5 \epsilon \epsilon' \pi \langle u_0, u_0 \rangle + \epsilon' \pi \langle u_0 - u , u_0 + u \rangle - \pi \langle u-u_0 , u -u_0 \rangle - c r_{\Lambda}^2(u).\]
Note that the norm of the first term is $5 \pi \epsilon \epsilon' |u_0|^2$.  The other terms are bounded above by $3 \pi \epsilon \epsilon' |u_0|$, $\pi \epsilon^2$, and $4 c \epsilon^2$.  In particular, since $u_0$ is an integral vector, and $\epsilon < \epsilon'$, the sum of the terms is bounded above by $10 \pi \epsilon \epsilon' |u_0|^2$ so long that $c$ is not too large.  It is also strictly positive since the first term dominates, so we conclude that the differential of our function does not attain the value $- 2 \pi u_0$ in the desired neighbourhood of $u_0$.

Since the norm of $dg$ is by assumption less than $\frac{2}{\epsilon}$, we conclude that the term corresponds to $dg$ has norm bounded by $20 \pi \epsilon' |u_0|^2$.  We can choose $\epsilon'$ small enough to ensure that this bound is much less than $1$, so that $dH_{c}$ is a bounded distance away from $- 2 \pi u_0$ as desired.
\end{proof}

The reason for the careful choices that we have made in the above paragraphs is as follows:  in order to interpolate between the differentials of two functions using cutoff functions without creating new critical points, we must know not only that the values of the differentials are close to each other, but also that the values of the functions are close.  This is exactly the reason why we need the comical choice of
\[ ( \epsilon' + 5 \epsilon \epsilon') \langle \psi(u), \psi(u) \rangle \]
as a term in our Hamiltonian.
\begin{cor} \label{prefloer1}
For sufficiently small $0<c$, $\cL_{c}(1)$ is an admissible Lagrangian that intersects the interior of $\cL$ in the set
\[ Q \cap \frac{1}{1+\epsilon'} \bZ^{n} \cong  Q \cap  \bZ^{n} .\]
Further, $\cL_{c}(1)$ agrees with a linear Lagrangian away from an $\epsilon$ neighbourhood of its boundary.
\end{cor}

If we want a model for $\cL(1)$ where the pair $(\cL, \cL(1))$ is negatively oriented, then we have to consider $\cL_{-c}(1)$.  However, following the above argument, we will find that the coefficient of $dr_{\Lambda}$ is negative in this case.  This means that we may be creating new intersection points if we try to interpolate between $\cL_{-c}(1)$ and the graph of $(1+\epsilon')\langle u, u \rangle$.  We can easily fix this by using  $(1-\epsilon')\langle u, u \rangle$ instead.  However,  this means that we must also use $- (\epsilon' + 5 \epsilon' \epsilon)$ in our expression for $H^{\partial}_c$.  The rest of the argument then carries through as desired and yields
\begin{cor} \label{prefloer2}
For sufficiently small $0<c$, $\cL_{-c}(1)$ is an admissible Lagrangian that intersects the interior of $\cL$ in the set
\[ (Q-\partial Q) \cap  \bZ^{n} .\]
Further, $\cL_{-c}(1)$ agrees with a linear Lagrangian away from an $\epsilon$ neighbourhood of its boundary.
\end{cor}

Since the Lagrangians $\cL_{\pm c}(1)$ are equivalent up to admissible Hamiltonian isotopy for any value of $c$ by \fullref{boundary_exists}, we revert to our old notation and refer to either of them as $\cL(1)$.

We can also use the same method to construct $\cL(j)$ for all $j$.  Recall that the first choice we made is that of a constant $0<\delta$, which allowed us to prove that no intersection point occurred in a neighbourhood of the boundary by bounding the norm of a certain vector by constants dominated by $\delta$.  It was sufficient, in order to construct $\cL(1)$, for $\delta$ to be much smaller than $1$, since we simply had to avoid any lattice points that are not in $Q$.  However, in order to construct $\cL(j)$, we will have to avoid all elements of the lattice $\frac{1}{j} \bZ^{n}$.  This means that we have to choose $\delta$ much smaller than $\frac{1}{j}$.  We can then use this to construct the Lagrangian boundary of $\cL(j)$ using the function $j H^{\partial}$, then carry through the rest of the steps of the preceding construction.  If $j$ in negative, we must work with $- \epsilon'$ as explained above. We conclude the folowing:
\begin{prop} \label{intersections}
For every integer $j \in \bZ$, there exists an admissible Lagrangian $\cL(j)$ such that the pair $(\cL,\cL(j))$ is positively oriented and their interior intersection points are in bijective correspondence with
\[ \frac{1}{j} \bZ^{N} \cap Q \]
if $j>0$ and with
\[ \frac{1}{j} \bZ^{N} \cap (Q- \partial Q) \]
if $j <0$.  Further, all interior intersection points occur on the complement of a small neighbourhood of the boundary. In this open set where the intersection points occur, $\cL(j)$ is given by the Hamiltonian
\[ - j \pi (1+ \epsilon') \langle u , u \rangle \]
for some small $\epsilon'$.
\end{prop}

\subsection{Computing Floer groups and products} \label{comp}
First we settle the issue of the Floer homology of $\cL$ with itself.
\begin{lem} \label{ordinary}
The only non-trivial Floer homology group of $\cL$ with itself is
\[ HF_{n}(\cL, \cL) = \bC .\]
\end{lem}
\begin{proof}
We take a small Hamiltonian isotopy and compute the Floer homology of $\cL$ with its image $\cL'$ under the isotopy.  Of course, this Hamiltonian isotopy must be chosen so that
$\cL'$ is admissible and the pair $(\cL, \cL')$ is positively
oriented.  To achieve this, we simply reproduce the arguments of the
previous section and interpolate between the graph of the differential of
\[ - \langle u , \eta u \rangle \]
with $\eta$ positive and the parallel transport of $\partial \cL$ along some curve $\gamma'$ which goes in the third quadrant. The point is that $\eta$
can be chosen so small that the differential of this function does not take integral values on $Q$. Further, the positivity of $\eta$ will guarantee that we can interpolate between the graph of this Hamiltonian and the parallel transport of the boundary in the direction of a ``positive'' curve $\gamma$ without creating any new intersection points with $\cL$.  The above quadratic function has a critical point of index $n$ at the
origin, so we reach the desired conclusion.  Note that this is
consistent with $H_{*}(D^{n}, S^{n-1})$.
\end{proof}

We can reinterpret \fullref{intersections} as follows:
\begin{lem}\label{floerhlgy}
The only non-vanishing Floer cohomology groups between $\cL$ and
$\cL(j)$ are given by
\begin{align*} HF^{0}(\cL,\cL(j)) &\cong \bigoplus_{p \in Q \cap \bZ^n}  \bC \cdot [p]\quad
\text{if $0 \leq j$\quad and,}\\
 HF^{n}(\cL,\cL(j)) &\cong \bigoplus_{p \in (Q - \partial Q) \cap \bZ^n} \bC \cdot [p] \quad
\text{if $j < 0$.}\end{align*}
\end{lem}
\begin{proof}
The fact that all intersections between $\cL$ and $\cL(j)$ occur in an open region where $\cL(j)$ is given by a linear Lagrangian allows us to compute the Floer cohomology as in \fullref{prelim}.
\end{proof}

We can now discuss the cup product.
\begin{lem}
Identifying the generators of Floer homology groups with lattice points as in \fullref{floerhlgy}, the product
\[
HF^{*}(\cL(l_1), \cL(l_2)) \otimes
HF^{*}(\cL(l_2), \cL(l_3)) \to
HF^{*}(\cL(l_1), \cL(l_3)) \]
is the same as the product for the Floer homology groups of the Lagrangians $\cL_{\infty}(l_i)$ which we used in \fullref{prelim}.  Concretely, for $l_1 \leq l_2 \leq l_3$, we have
 \[ [ p ] \otimes [ q ] \mapsto \left[ \frac{(l_2-l_1) p + (l_3-l_2) q}{l_3-l_1} \right]\]
as in \fullref{prelim_product}.
\end{lem}
\begin{proof}
First, we observe that we can choose the constant $\delta$ is the previous section so that our models for $\cL(l_i)$ have all their interior intersection points in the open region where these Lagrangians are linear.  This means that all the computations of \fullref{prelim} hold.  We now justify the perturbation arguments used in that section.

It should be clear now that the Maslov index computation for triple intersections that we performed is legitimate, even when the intersection point occurs at the boundary of $\cL_{\infty}(l_i)$, since in $\cL(l_i)$ these are now interior intersection points, and can be perturbed to achieve transversality without affecting the admissibility condition.  It remains to address the situation where two of the three Lagrangians are equal. The most important thing is to check that our generator of $HF^*(\cL, \cL)$ acts as the identity. Say we are trying to compute
\begin{equation*} HF^{*}( \cL , \cL) \otimes HF^{*}( \cL , \cL(j)) \to HF^{*}( \cL , \cL(j)) . \end{equation*}
Now choose an intersection point $[q]$ between $\cL$ and $\cL(j)$, and perturb $\cL$ to a Lagrangian $\cL'$ as in \fullref{ordinary}. Note that we can ensure that the intersection point $p$ between $\cL$ and $\cL'$ occurs near $q$.  We claim that the cup product is given by
\[ [p] \otimes [q'] \to [q] \]
where $[q']$ is the corresponding nearby intersection point between $\cL'$ and $\cL(j)$.  In effect, we have reduced everything to a Maslov index computation.  The point is that we can choose all these intersections to occur at the same point, and $\cL'$ to have slope $-\eta$.  For sufficiently small $\eta$, the numbers $( 0 , \eta, j)$ satisfy Conditions \eqref{cup1} or \eqref{cup2} regardless of what the sign of $j$ is.  The usual obstruction in the first homology of the torus fibre allows us to conclude that there can be no other holomorphic triangle.
\end{proof}

\subsection{Completing the argument} \label{end}

Let $X$ be the toric variety of the Introduction with an ample line bundle $\cO(1)$. The last section,  together with Equation \eqref{proj}, allows us to conclude that if $j \leq l$, then
\[ \Hom_{*} ( \cO(j) , \cO(l) ) \cong HF^{*} (\cL(j), \cL(l) ), \]
and that this isomorphism is compatible with cup product. This proves \fullref{main}.  

We now use Serre duality to compute the case $j > l$.  To simplify the notation, we let $j=0$, so $l<0$ and 
\[ \Hom_{*} ( \cO , \cO(l) ) \cong  (\Hom_{n-*} ( \cO , \cO(-l) \otimes \kappa)) \check{\,} . \]
Where $\kappa$ is the canonical bundle.  Since the piecewise linear function which defines $\kappa$ takes the value $-1$ on each primitive vertex of a $1$--cone of $\Delta$, \acite{fulton}*{Section 4.3}, we see that our requirement that the interior of $Q$ have an integral point (See \fullref{int}), is equivalent to the fact that the piecewise linear function defining $ \cO(-l) \otimes \kappa$ is convex.  Its sections are in fact given by integral points which satisfy
\[ \langle v_{i}, y \rangle \leq \phi(v_i)-1 \]
for every primitive vertex $v_i$ of the $1$--cones of $\Delta$.  But these are exactly the interior lattice points of $-lQ$. We have therefore proved that
\[ \Hom_{*} ( \cO(j) , \cO(l) ) \cong HF^{*} (\cL(j), \cL(l) ) \]
for all values of $j$ and $l$.  One may check that this isomorphism is indeed compatible with the cup product (which, for negative powers is interpreted through Serre duality).  We check this for the case $l>0$, $j < -l$ by computing
\[ \Hom_{*} ( \cO , \cO(l) ) \otimes  \Hom_{*} ( \cO , \cO(j) ) \to \Hom_{*} ( \cO , \cO(l+j) ) .\]
Note that Serre duality suggests that we should reduce this computation to that of
\[
\Hom_{*} ( \cO , \cO(l) ) \otimes  \Hom_{*} ( \cO , \cO(-l-j) \otimes \kappa )^{\check{\,}} \to \Hom_{*} ( \cO , \cO(-j) \otimes \kappa), \]
in which all $\Hom$'s are concentrated in degree $0$.  The formula for the product is
\[ [p] \otimes [q] \to [p+q] ,\]
which then dualizes to
\[ [p] \otimes [r] \check{\,} = \sum_{q \in (l+j) (Q - \partial Q) \cap \bZ^{n} \ni  p+q=r } [q] \check{\,} , \]
where the summation in the right hand side is either empty or consists of a unique term.  The usual change of perspective from lattice points to $(l+j) Q$ to $\frac{1}{l+j}$ lattice points of $Q$ yields the isomorphism with the product on Floer cohomology.

\bibliographystyle{gtart}
\bibliography{link}

\begin{thebibliography}{}
\providecommand\bibmarginpar{\leavevmode\marginpar}
\def\urlstyle#1{{\tt #1}}

\bibitem{beilinson}
\textbf{A\,A Be{\u\i}linson}, \emph{Coherent sheaves on $\mathbf{P}^{n}$ and
  problems in linear algebra}, Funktsional. Anal. i Prilozhen. 12 (1978) 68--69
  \xox{MR}{509388}

\bibitem{deSilva}
\textbf{V De~Silva}, \emph{Products in the symplectic Floer homology of
  Lagrangian intersections}, PhD thesis, Oxford University (1998)

\bibitem{donaldson}
\textbf{S\,K Donaldson}, \emph{Symplectic submanifolds and almost-complex
  geometry}, J. Differential Geom. 44 (1996) 666--705 \xox{MR}{1438190}

\bibitem{EKL}
\textbf{M Einsiedler}, \textbf{M Kapranov}, \textbf{D Lind},
  \emph{Non-archimedean amoebas and tropical varieties}
  \xox{arXiv}{math.AG/0408311}

\bibitem{eliashberg-gromov}
\textbf{Y Eliashberg}, \textbf{M Gromov}, \emph{Convex symplectic manifolds},
  from: ``Several complex variables and complex geometry, Part 2 (Santa Cruz,
  CA, 1989)'', Proc. Sympos. Pure Math. 52, Amer. Math. Soc., Providence, RI
  (1991)  135--162 \xox{MR}{1128541}

\bibitem{floer}
\textbf{A Floer}, \emph{Morse theory for Lagrangian intersections}, J.
  Differential Geom. 28 (1988) 513--547 \xox{MR}{965228}

\bibitem{FHS}
\textbf{A Floer}, \textbf{H Hofer}, \textbf{D Salamon}, \emph{Transversality in
  elliptic Morse theory for the symplectic action}, Duke Math. J. 80 (1995)
  251--292 \xox{MR}{1360618}

\bibitem{FO}
\textbf{K Fukaya}, \textbf{Y-G Oh}, \emph{Zero-loop open strings in the
  cotangent bundle and {M}orse homotopy}, Asian J. Math. 1 (1997) 96--180
  \xox{MR}{1480992}

\bibitem{FOOO3}
\textbf{K Fukaya}, \textbf{Y-G Oh}, \textbf{H Ohta}, \textbf{K Ono},
  \emph{Lagrangian intersection Floer theory - anomaly and obstruction},
  unpublished manuscript (2000)
\ Available at \setbox0\hbox{\makeatletter\@url
{http://www.math.kyoto-u.ac.jp/~fukaya/}}
\href{http://www.math.kyoto-u.ac.jp/~fukaya/}
{\unhbox0}

\bibitem{fulton}
\textbf{W Fulton}, \emph{Introduction to toric varieties}, Annals of
  Mathematics Studies 131, Princeton University Press, Princeton, NJ (1993)
  \xox{MR}{1234037}\ \ The William H. Roever Lectures in Geometry

\bibitem{GKZ}
\textbf{I\,M Gel'fand}, \textbf{M\,M Kapranov}, \textbf{A\,V Zelevinsky},
  \emph{Discriminants, resultants, and multidimensional determinants},
  Mathematics: Theory \& Applications, Birkh\"auser, Boston (1994)
  \xox{MR}{1264417}

\bibitem{kont-ENS}
\textbf{M Kontsevich}, \emph{Lectures at ENS Paris, Spring 1998}\ \ Notes by J
  Bellaiche, J-F Dat, I Marin, G Racinet and H Randriambololona

\bibitem{kontsevich}
\textbf{M Kontsevich}, \emph{Homological algebra of mirror symmetry}, from:
  ``Proceedings of the International Congress of Mathematicians, Vol.\ 1, 2
  (Z\"urich, 1994)'', Birkh\"auser, Basel (1995)  120--139 \xox{MR}{1403918}

\bibitem{leung}
\textbf{N\,C Leung}, \emph{Mirror symmetry without corrections}
  \xox{arXiv}{math.DG/0009235}

\bibitem{mikhalkin}
\textbf{G Mikhalkin}, \emph{Decomposition into pairs-of-pants for complex
  algebraic hypersurfaces}, Topology 43 (2004) 1035--1065 \xox{MR}{2079993}

\bibitem{ruan}
\textbf{W-D Ruan},
  \href{http://projecteuclid.org/getRecord?id=euclid.jsg/1092403030}
  {\emph{Lagrangian torus fibration of quintic {C}alabi--{Y}au hypersurfaces
  II: {T}echnical results on gradient flow construction}}, J. Symplectic Geom.
  1 (2002) 435--521 \xox{MR}{1959057}

\bibitem{rullgaard}
\textbf{H Rullg{\aa}rd}, \emph{Polynomial amoebas and convexity}, preprint,
  Stockholm University  (2001)

\bibitem{seidelGL}
\textbf{P Seidel}, \emph{Graded {L}agrangian submanifolds}, Bull. Soc. Math.
  France 128 (2000) 103--149 \xox{MR}{1765826}

\end{thebibliography}

\end{document}